\documentclass[11pt]{article}
\usepackage{amsthm,amstext,amssymb,amsmath,graphicx,amsfonts,caption,subcaption,multirow,epigraph,mathrsfs}
\usepackage[all]{xy}
\usepackage[utf8]{inputenc}
\usepackage{booktabs}
\usepackage{tikz}
\usetikzlibrary{decorations.pathreplacing}
\usetikzlibrary{arrows.meta, positioning, shapes.geometric, calc}
\usetikzlibrary{matrix}
\usetikzlibrary{arrows,topaths}
\usepackage{placeins}
\usepackage{mathrsfs}
\usepackage{multicol}
\usepackage[numbers]{natbib}
\usepackage[pagebackref=false,colorlinks,linkcolor=blue,citecolor=blue]{hyperref}
\usepackage{tikz-cd}
\usepackage{xcolor,soul}
\theoremstyle{plain}
\newtheorem{theorem}{Theorem}[section]
\newtheorem{lemma}[theorem]{Lemma}
\newtheorem{proposition}[theorem]{Proposition}

\theoremstyle{definition}
\newtheorem{definition}[theorem]{Definition}

\newtheorem{example}[theorem]{Example}
\newtheorem{examples}[theorem]{Examples}

\theoremstyle{remark}
\newtheorem{remark}{\sc Remark}
\makeatother
\makeatletter
\def\namedlabel#1#2{\begingroup
   \def\@currentlabel{#2}%
   \label{#1}\endgroup}
\makeatother
\date{}
\title{\bf  Sheaf of residuated lattices}\vspace{.25 in}
\author{ \vspace{.25 in} {\bf Saeed Rasouli}\\
Department of Mathematics, Persian Gulf University, \\Bushehr, Iran \\
{\tt srasouli@pgu.ac.ir }\\}

\begin{document}
 \maketitle
 \begin{abstract}
 This paper explores the interface between algebra, topology, and logic by developing the theory of sheaves and \'{e}tal\'{e} spaces for \emph{residuated lattices}—algebraic structures central to substructural and fuzzy logics. We construct stalkwise-residuated \'{e}tal\'{e} spaces and demonstrate that they form a subcategory of the category of \'{e}tal\'{e} spaces of sets. A categorical and topological characterization of the sheaf condition is presented, with particular emphasis on filters, congruences, and the topologies induced on prime spectra.

We develop presheaves and sheaves valued in the category of residuated lattices, define sheafification via direct limits, and establish an equivalence between sections of \'{e}tal\'{e} spaces and sheaf models. Illustrative examples—including the {\L}ukasiewicz structure and the Sierpi\'nski space—demonstrate the modeling of local truth in both spatial and logical contexts.

This paper concludes by establishing a functorial correspondence between algebraic \'{e}tal\'{e} spaces and logical sheaves, thereby providing a categorical framework for context-sensitive reasoning. Specifically, we demonstrate that the category of sheaves of residuated lattices over a topological space~$\mathscr{B}$ is \emph{reflective} within the category of presheaves of residuated lattices over~$\mathscr{B}$. Furthermore, we show that this category is \emph{equivalent} to the category of \'{e}tal\'{e} spaces of residuated lattices over~$\mathscr{B}$.\footnote{2020 Mathematics Subject Classification: 06F99,55N30  \\
{\it Key words and phrases}: sheaf of residuated lattices; presheaf of residuated lattices;\'{e}tal\'{e} space of residuated lattice.}
\end{abstract}


\section{Introduction}

Sheaf theory is a powerful framework in modern mathematics, offering a systematic way to relate local data to global structures. Its origins can be traced to the work of Jean Leray during World War II, when he was held as a prisoner of war in the German camp Oflag XVII. During this period, Leray developed foundational ideas while investigating solutions to differential equations on manifolds, leading to the introduction of sheaves as a means of organizing local solutions and understanding their global behavior \cite{leray1945forme, leray1946lanneau, leray1949homologie}. His pioneering insights laid the groundwork for what would become one of the most influential theories in algebraic topology, differential geometry, and later algebraic geometry. Leray's ideas were later expanded and refined in the influential S\'{e}minaire Henri Cartan (1948–1964), led by Henri Cartan in collaboration with eminent mathematicians such as Jean-Pierre Serre, Michel Lazard, and Roger Godement \citep{cartan1950faisceaux}.

A foundational development in this seminar was the notion of an \emph{\'etal\'e space} (\emph{espace \'etal\'e}), closely tied to the concept of a sheaf. It is widely acknowledged that Henri Cartan, often in collaboration with Michel Lazard, introduced the geometric notion of a sheaf as an \'etale space in \emph{Expos\'e 14} of the seminar, under the term \emph{faisceau} \citep{cartan1950faisceaux}. According to \citet[p.~7]{Gray1979}, Lazard’s faisceau on a regular topological space $X$ consists of a map $p: F \to X$ satisfying:
\begin{enumerate}
  \item For each $x \in X$, the fiber $p^{-1}(x) = F_x$ is a $K$-module;
  \item The total space $F$ carries a (possibly non-Hausdorff) topology such that the algebraic operations on $F$ are continuous, and $p$ is a local homeomorphism.
\end{enumerate}

The term \emph{\'{e}tal\'{e} space} was later coined and popularized by Roger Godement, a central participant in the Cartan seminar. In his influential treatise on general topology \cite{godement1958topologie}, Godement offered a systematic and elegant treatment of \'etale spaces, contributing significantly to their formalization. His exposition helped shape the understanding of sheaf-theoretic structures in algebraic topology and algebraic geometry.

The formal notion of a \emph{presheaf} was introduced by \citet{borelcohomologie}, complementing the development of sheaf theory by allowing the description of local sections over open subsets without necessarily satisfying gluing conditions. This distinction proved essential in the broader categorical framework that emerged.

Over the decades, sheaf theory and the associated concept of \'etal\'e spaces have found applications far beyond their origins in topology and algebraic geometry. Contributions by mathematicians such as \citet{Tennison1976,Bredon1997,grothendieck1960elements,maclane2012sheaves,brown2020sheaf,rosiak2022sheaf} have deepened and expanded the theory in foundational and categorical directions. More recently, \citet{hohle2007fuzzy} demonstrated that fuzzy set theory—concerned with modeling vague or uncertain information—can be naturally described within a sheaf-theoretic framework. This development underscores the far-reaching impact of sheaf theory, illustrating how ideas originating in pure mathematics continue to inform and shape applications in logic, information theory, and beyond. For more historical notes about sheaf spaces, the interested reader can be referred to \citet{Gray1979,fasanelli1981creation,miller2000leray}.


Classical logic serves as a foundational framework for information processing, particularly for reasoning based on definite, crisp information. However, many real-world contexts involve uncertainty, vagueness, or gradations of truth. To address these phenomena, it is necessary to develop logical systems that extend beyond classical logic. As a result, a wide array of non-classical logics has been introduced and extensively studied. These logics now play a central role in fields such as computer science and artificial intelligence, where managing uncertainty and imprecision is essential.

In parallel with this logical development, various algebraic structures have emerged as semantic models for non-classical logics. Among these are residuated lattices, divisible residuated lattices, MTL-algebras, Girard monoids, BL-algebras, MV-algebras, and G\"odel algebras. Residuated lattices, in particular, form a versatile class of algebras that model many substructural and fuzzy logics, providing a coherent algebraic semantics for systems lacking certain structural rules. Their flexibility and rich structure have made them a cornerstone in the algebraic study of non-classical reasoning \citep{galatos2007residuated, hajek2013metamathematics, turunen1999mathematics}.

One especially notable application of residuated lattices lies in \emph{fuzzy logic}, where truth values are drawn from a continuum rather than from binary values. In this setting, residuated lattices provide the semantic underpinning for key operations such as conjunction, implication, and aggregation. The adjointness property of the monoidal operation and residuation is critical in preserving logical coherence when reasoning with partial truth. Systems such as \emph{MTL} (Monoidal t-norm based Logic) and \emph{BL} (Basic Logic) rely explicitly on the algebraic properties of residuated lattices to model graded entailment and approximate reasoning. This makes residuated lattices indispensable for capturing the behavior of fuzzy connectives and t-norm-based semantics.

In Gentzen-style sequent calculi, a \emph{structural rule} is one that governs the manipulation of sequents independently of logical connectives. Substructural logics are defined as systems where one or more of the following structural rules are absent:
\begin{itemize}
  \item \textbf{Weakening:}
  \[
  \frac{\Gamma, \Delta \Rightarrow \varphi}{\Gamma, \alpha, \Delta \Rightarrow \varphi}
  \]
  \item \textbf{Contraction:}
  \[
  \frac{\Gamma, \alpha, \alpha, \Delta \Rightarrow \varphi}{\Gamma, \alpha, \Delta \Rightarrow \varphi}
  \]
  \item \textbf{Exchange:}
  \[
  \frac{\Gamma, \alpha, \beta, \Delta \Rightarrow \varphi}{\Gamma, \beta, \alpha, \Delta \Rightarrow \varphi}
  \]
\end{itemize}

In this paper, we focus on \emph{commutative residuated lattices}, which arise naturally as the algebraic semantics of substructural logics lacking the contraction rule.

Although residuated lattices are now prominent in logic, their origins are algebraic rather than logical. The foundational work on such structures dates back to \citet{krull1924axiomatische}, who studied the decomposition of ideals in rings. Subsequent developments by M.~Ward and R.~P.~Dilworth in a series of influential papers \citep{dilworth1938abstract, dilworth1939non, ward1937residuation, ward1938structure, ward1940residuated, ward1938residuated, ward1939residuated} established residuated lattices as fundamental tools in the abstract theory of ideal lattices. Over time, these lattices have been studied under various names, including BCK-lattices \citep{hohle1995commutative}, full BCK-algebras \citep{krull1924axiomatische}, $FL_{ew}$-algebras \citep{okada1999finite}, and integral, residuated, commutative $\ell$-monoids \citep{blok1989algebraizable}.


Given the utility of sheaf theory for relating local and global structures, and the importance of residuated lattices in modeling uncertainty, it is natural to explore the intersection of these two frameworks. This leads to the notion of a \emph{sheaf of residuated lattices}, an analogue of sheaves of rings that enables a topological understanding of logical systems. Such sheaves are particularly relevant in fuzzy logic, where one often models locally varying degrees of truth.

The motivation for defining a sheaf of residuated lattices stems from the desire to capture both the \emph{local variation of truth values}---as encountered in fuzzy environments---and the \emph{global logical coherence} that emerges from the interrelation of these local structures. In many logical systems, particularly those addressing vagueness or partial knowledge, different regions or contexts may exhibit distinct inferential behavior or differing algebraic semantics. A sheaf-theoretic approach allows one to formalize this locality in a rigorous way, enabling the patching together of locally valid logics into a globally meaningful whole.

Moreover, such a framework facilitates the study of logical semantics that are not uniform across a space but instead vary continuously or discretely with respect to a topology. This is crucial in modeling contexts such as sensor networks, distributed information systems, or spatial reasoning in artificial intelligence---where local truths must be synthesized coherently. In fuzzy logic, where truth values are typically drawn from a residuated lattice, the sheaf perspective allows for the systematic integration of variable or context-dependent truth degrees, thus enriching the expressiveness and applicability of the logic.

From a practical standpoint, sheaves of residuated lattices have promising applications in real-world domains where knowledge is inherently local and uncertain. For example, in \emph{geospatial information systems}, truth values about terrain properties, land use, or sensor outputs often vary across regions and depend on local data quality. Modeling such environments using sheaves of residuated lattices provides a principled way to aggregate and reconcile these localized truth values while preserving logical consistency. Similarly, in \emph{context-aware computing}---such as smart environments or adaptive user interfaces---logical decisions must account for varying, imprecise contextual inputs. The sheaf-based structure enables the seamless integration of these contextual inferences into a coherent global state, enabling robust and interpretable reasoning in dynamic, uncertain settings.

Related ideas have appeared across various domains: sheaves of lattices \citep{brezuleanu1969duale}, sheaves of universal algebras \citep{davey1973sheaf}, sheaves of MV-algebras \citep{filipoiu1995compact}, sheaves of BL-algebras \citep{di2003compact}, and sheaves of Banach spaces \citep{hofmann2006sheaf}. These studies underscore the flexibility and broad applicability of sheaf-theoretic methods.

In this paper, we initiate a systematic study of sheaves of residuated lattices. Our goal is to uncover categorical and topological structures that generalize known results for rings, while adopting a more direct approach that avoids reliance on the full machinery of Grothendieck-style sheaf theory. We aim to clarify how uncertainty and locality interact in logical systems by using algebraic and topological tools in tandem.


This paper is organized into seven sections as follows. 

In Section~\ref{sec2}, we review the foundational concepts relevant to our development. The first subsection covers essential definitions, properties, and results concerning residuated lattices, while the second subsection revisits the theory of \'etale spaces of sets. These preliminaries serve as the groundwork for the subsequent sections.

Although many of the results in the second subsection are well known and appear in standard references on sheaf theory, we have chosen to include full proofs for the sake of completeness and to ensure self-containment. Moreover, the presentation and order of propositions and the style of proofs have been tailored to align with the general framework and methodological principles adopted in this work.

Section~\ref{sec3} is devoted to the study of \'{e}tal\'e spaces of residuated lattices over a fixed topological space. We introduce the notion of morphisms between such structures and investigate their categorical properties. In particular, we show that for a given space~$\mathscr{B}$, the class of \'{e}tal\'e spaces of residuated lattices over~$\mathscr{B}$, together with their morphisms, forms a subcategory of the category of \'{e}tal\'{e} spaces of sets over \(\mathscr{B}\) (Theorem~\ref{etressubcatetsp}).

Section~\ref{sec4} focuses on presheaves and sheaves of residuated lattices. We begin by introducing the notions of a presheaf and a sheaf of residuated lattices over a topological space, accompanied by illustrative examples that clarify these definitions. The section culminates in the construction of a functor from the category of \'{e}tal\'e spaces of residuated lattices to the category of sheaves of residuated lattices (Theorem~\ref{etalesheafunc}). Notably, examples such as the {\L}ukasiewicz structure and the Sierpi\'nski space are examined to illustrate how these frameworks model local truth in both spatial and logical contexts.

Section~\ref{sec5} presents the sheafification process via direct limits, which transforms a presheaf of residuated lattices into a sheaf of residuated lattices. We begin by decomposing a presheaf into its germs (Theorem~\ref{shefifipro1}), which capture the finest local behavior of sections at individual points. Based on this decomposition, we construct a functor from the category of presheaves of residuated lattices over a topological space~$\mathscr{B}$ to the category of \'{e}tal\'e spaces of residuated lattices over~$\mathscr{B}$ (Theorem~\ref{preshetalfunc}).

A sheaf is then reconstructed by assembling these germs into a globally coherent structure, subject to the condition that local sections agree on overlaps. The resulting sheaf consists precisely of those sections that can be consistently patched together from compatible local data.

Section~\ref{sec6} investigates the categorical relationships among \'{e}tal\'e spaces, presheaves, and sheaves of residuated lattices over a fixed topological space. We prove that for a given space~$\mathscr{B}$, the category of sheaves of residuated lattices over~$\mathscr{B}$ is reflective in the category of presheaves of residuated lattices over~$\mathscr{B}$ (Theorem~\ref{reflshbpshb}). Furthermore, we establish that the category of \'{e}tal\'e spaces of residuated lattices over~$\mathscr{B}$ is equivalent to the category of sheaves of residuated lattices over~$\mathscr{B}$ (Theorem~\ref{eqetshtheo}).

Finally, Section~\ref{sec7} presents the conclusion of the paper, along with a logical interpretation of the main results. The section also outlines several open problems and directions for future research, emphasizing the potential for further exploration at the intersection of algebra, logic, and topology.

\section{Preliminaries}\label{sec2}
%
%
\subsection{residuated lattices}

Throughout this paper, all residuated lattices are assumed to be integral, commutative, and non-degenerate unless stated otherwise. In this section, we recall some basic definitions and properties that will be used in subsequent sections.

An algebra \(\mathfrak{A} = (A; \vee, \wedge, \odot, \rightarrow, 0, 1)\) is called a \textit{residuated lattice} if the following conditions are satisfied:
\begin{itemize}
	\item \(\ell(\mathfrak{A}) = (A; \vee, \wedge, 0, 1)\) is a bounded lattice;
	\item \((A; \odot, 1)\) is a commutative monoid;
	\item \((\odot, \rightarrow)\) forms an adjoint pair.
\end{itemize}
On any residuated lattice \(\mathfrak{A}\), a unary operation \(\neg\) can be defined as \(\neg a := a \rightarrow 0\). The class of residuated lattices is equational and thus forms a variety. For an extensive survey of residuated lattices, we refer the reader to \cite{galatos2007residuated}.

\begin{example}\label{exa4}
Consider the lattice $A_4 = \{0, a, b, 1\}$ with the Hasse diagram shown in Figure \ref{figa4}. Routine calculations reveal that $\mathfrak{A}_4 = (A_4; \vee, \wedge, \odot, \rightarrow, 0, 1)$ is a residuated lattice. The commutative operation $\odot$ is defined by Table \ref{taba4}, and the operation $\rightarrow$ is articulated by $x \rightarrow y = \bigvee \{z \in A_4 \mid x \odot z \leq y\}$, for any $x, y \in A_4$.
\FloatBarrier
\begin{table}[ht]
\begin{minipage}[b]{0.56\linewidth}
\centering
\begin{tabular}{ccccccc}
\hline
$\odot$ & 0 & a & b & 1 \\ \hline
0       & 0 & 0 & 0 & 0 \\
        & a & a & 0 & a \\
        &   & b & b & b \\
        &   &   & 1 & 1 \\ \hline
\end{tabular}
\caption{Cayley table for ``$\odot$" of $\mathfrak{A}_4$}
\label{taba4}
\end{minipage}\hfill
\begin{minipage}[b]{0.6\linewidth}
\centering
  \begin{tikzpicture}[>=stealth',semithick,auto]
    \tikzstyle{subj} = [circle, minimum width=6pt, fill, inner sep=0pt]
    \tikzstyle{obj}  = [circle, minimum width=6pt, draw, inner sep=0pt]

    \tikzstyle{every label}=[font=\bfseries]

    \node[subj,  label=below:0] (0) at (0,0) {};
    \node[subj,  label=below:a] (a) at (-1,1) {};
    \node[subj,  label=below:b] (b) at (1,1) {};
    \node[subj,  label=right:1] (1) at (0,2) {};

    \path[-]   (0)    edge                node{}      (a);
    \path[-]   (0)    edge                node{}      (b);
    \path[-]   (a)    edge                node{}      (1);
    \path[-]   (b)    edge                node{}      (1);
\end{tikzpicture}
\captionof{figure}{Hasse diagram of $\mathfrak{A}_{4}$}
\label{figa4}
\end{minipage}
\end{table}
\FloatBarrier
\end{example}
\begin{example}\label{exa6}
Consider the lattice $A_6=\{0,a,b,c,d,1\}$ with the Hasse diagram shown in Figure \ref{figa6}. Routine calculations reveal that $\mathfrak{A}_6=(A_6;\vee,\wedge,\odot,\rightarrow,0,1)$ is a residuated lattice. The commutative operation $\odot$ is defined by Table \ref{taba6}, and the operation $\rightarrow$ is articulated by $x \rightarrow y = \bigvee \{z \in A_6 \mid x \odot z \leq y\}$, for any $x, y \in A_6$.
\FloatBarrier
\begin{table}[ht]
\begin{minipage}[b]{0.56\linewidth}
\centering
\begin{tabular}{ccccccc}
\hline
$\odot$ & 0 & a & b & c & d & 1 \\ \hline
0       & 0 & 0 & 0 & 0 & 0 & 0 \\
        & a & a & a & 0 & a & a \\
        &   & b & a & 0 & a & b \\
        &   &   & c & c & c & c \\
        &   &   &   & d & d & d \\
        &   &   &   &   & 1 & 1 \\ \hline
\end{tabular}
\caption{Cayley table for ``$\odot$" of $\mathfrak{A}_6$}
\label{taba6}
\end{minipage}\hfill
\begin{minipage}[b]{0.6\linewidth}
\centering
  \begin{tikzpicture}[>=stealth',semithick,auto]
    \tikzstyle{subj} = [circle, minimum width=6pt, fill, inner sep=0pt]
    \tikzstyle{obj}  = [circle, minimum width=6pt, draw, inner sep=0pt]

    \tikzstyle{every label}=[font=\bfseries]

    \node[subj,  label=below:0] (0) at (0,0) {};
    \node[subj,  label=below:c] (c) at (-1,1) {};
    \node[subj,  label=below:a] (a) at (1,.5) {};
    \node[subj,  label=below right:b] (b) at (1,1.5) {};
    \node[subj,  label=below:d] (d) at (0,2) {};
    \node[subj,  label=below right:1] (1) at (0,3) {};

    \path[-]   (0)    edge                node{}      (a);
    \path[-]   (a)    edge                node{}      (b);
    \path[-]   (0)    edge                node{}      (c);
    \path[-]   (c)    edge                node{}      (d);
    \path[-]   (b)    edge                node{}      (d);
    \path[-]   (d)    edge                node{}      (1);
\end{tikzpicture}
\captionof{figure}{Hasse diagram of $\mathfrak{A}_{6}$}
\label{figa6}
\end{minipage}
\end{table}
\FloatBarrier
\end{example}
\begin{example}\label{exa8}
Consider the lattice $A_8=\{0,a,b,c,d,e,f,1\}$ with the Hasse diagram shown in Figure \ref{figa8}. Routine calculations reveal that $\mathfrak{A}_8=(A_8;\vee,\wedge,\odot,\rightarrow,0,1)$ is a residuated lattice. The commutative operation $\odot$ is defined by Table \ref{taba8}, and the operation $\rightarrow$ is articulated by $x \rightarrow y = \bigvee \{z \in A_8 \mid x \odot z \leq y\}$, for any $x, y \in A_8$.
\FloatBarrier
\begin{table}[ht]
\begin{minipage}[b]{0.56\linewidth}
\centering
\begin{tabular}{ccccccccc}
\hline
$\odot$ & 0 & a & b & c & d & e & f & 1 \\ \hline
0       & 0 & 0 & 0 & 0 & 0 & 0 & 0 & 0 \\
        & a & 0 & a & a & a & a & a & a \\
        &   & b & 0 & 0 & 0 & 0 & b & b \\
        &   &   & c & c & a & c & a & c \\
        &   &   &   & d & a & a & d & d \\
        &   &   &   &   & e & c & d & e \\
        &   &   &   &   &   & f & f & f \\
        &   &   &   &   &   &   & 1 & 1 \\ \hline
\end{tabular}
\caption{Cayley table for ``$\odot$" of $\mathfrak{A}_8$}
\label{taba8}
\end{minipage}\hfill
\begin{minipage}[b]{0.6\linewidth}
\centering
  \begin{tikzpicture}[>=stealth',semithick,auto]
    \tikzstyle{subj} = [circle, minimum width=6pt, fill, inner sep=0pt]
    \tikzstyle{obj}  = [circle, minimum width=6pt, draw, inner sep=0pt]

    \tikzstyle{every label}=[font=\bfseries]

    \node[subj,  label=below:0] (0) at (0,0) {};
    \node[subj,  label=below:a] (a) at (-1,1) {};
    \node[subj,  label=below:b] (b) at (1,1) {};
    \node[subj,  label=below:c] (c) at (-2,2) {};
    \node[subj,  label=below:d] (d) at (0,2) {};
    \node[subj,  label=below:e] (e) at (-1,3) {};
    \node[subj,  label=below:f] (f) at (1,3) {};
    \node[subj,  label=below:1] (1) at (0,4) {};

    \path[-]   (0)    edge                node{}      (a);
    \path[-]   (0)    edge                node{}      (b);
    \path[-]   (b)    edge                node{}      (d);
    \path[-]   (d)    edge                node{}      (f);
    \path[-]   (f)    edge                node{}      (1);
    \path[-]   (a)    edge                node{}      (d);
    \path[-]   (a)    edge                node{}      (c);
    \path[-]   (c)    edge                node{}      (e);
    \path[-]   (d)    edge                node{}      (e);
    \path[-]   (e)    edge                node{}      (1);
\end{tikzpicture}

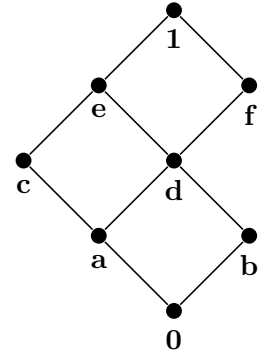
\captionof{figure}{Hasse diagram of $\mathfrak{A}_{8}$}
\label{figa8}
\end{minipage}
\end{table}
\FloatBarrier
\end{example}
\begin{example}\label{lukas}
Let $n$ be a fixed natural number. As noted by \citet[p. 11, Example 2]{turunen1999mathematics}, the structure $\L = ([0,1]; \max, \min, \odot, \rightarrow, 0, 1)$ forms a residuated lattice, referred to as the \textit{generalized {\L}ukasiewicz structure}. In this structure, the commutative operation $\odot$ is defined as 
\[
x \odot y = \left(\max\{0, x^n + y^n - 1\}\right)^{\frac{1}{n}},
\]
and the implication operation $\rightarrow$ is given by 
\[
x \rightarrow y = \min\left\{1, \left(1 - x^n + y^n\right)^{\frac{1}{n}}\right\},
\]
for all $x, y \in [0,1]$. 

This structure serves as a foundational model in {\L}ukasiewicz fuzzy logic, where the operations $\odot$ and $\rightarrow$ represent conjunction and implication, respectively, capturing the behavior of truth values in many-valued logical systems.

\end{example}

Let \(\mathfrak{A}\) be a residuated lattice. A non-empty subset \(F\) of \(A\) is called a \textit{filter} of \(\mathfrak{A}\) if the following conditions hold:
\begin{itemize}
	\item For all \(x, y \in F\), \(x \odot y \in F\);
	\item For any \(x \in F\) and \(y \in A\), \(x \vee y \in F\).
\end{itemize}
The set of all filters of \(\mathfrak{A}\) is denoted by \(\mathscr{F}(\mathfrak{A})\). $\{1\}$ and $A$ are called trivial filters. A filter \(F\) of \(\mathfrak{A}\) is called \textit{proper} if \(F \neq A\).

For any subset \(X\) of \(A\), the \textit{filter of \(\mathfrak{A}\) generated by \(X\)} is denoted by \(\mathscr{F}^{\mathfrak{A}}(X)\); here, \(\mathfrak{A}\) may be omitted when understood from context. For each \(x \in A\), the filter generated by \(\{x\}\) is denoted by \(\mathscr{F}(x)\) and is called a \textit{principal filter}. The set of all principal filters is denoted by \(\mathscr{PF}(\mathfrak{A})\).

Following \citet[\S 5.7]{gratzer2011lattice}, a join-complete lattice \(\mathfrak{A}\) is called a \textit{frame} if it satisfies the join infinite distributive law (JID), i.e., for any \(a \in A\) and \(S \subseteq A\),
\[
a \wedge \bigvee S = \bigvee \{a \wedge s \mid s \in S\}.
\]
A frame \(\mathfrak{A}\) is called \textit{complete} if \(\mathfrak{A}\) is a complete lattice.

According to \cite{galatos2007residuated}, \((\mathscr{F}(\mathfrak{A}); \cap, \veebar, \textbf{1}, A)\) is a complete frame, where for any \(\mathcal{F} \subseteq \mathscr{F}(\mathfrak{A})\), \(\veebar \mathcal{F} = \mathscr{F}(\cup \mathcal{F})\).

\begin{example}\label{filterexa}
Consider the residuated lattice $\mathfrak{A}_4$ from Example \ref{exa4}, the residuated lattice $\mathfrak{A}_6$ from Example \ref{exa6}, and the residuated lattice $\mathfrak{A}_8$ from Example \ref{exa8}. The sets of their filters are presented in Table \ref{tafiex}.
\begin{table}[h]
\centering
\begin{tabular}{ccl}
\hline
                 & \multicolumn{2}{c}{Filters}                                       \\ \hline
$\mathfrak{A}_4$ & \multicolumn{2}{c}{$F_{1}=\{1\},F_{2}=\{a,1\},F_{3}=\{b,1\},F_{4}=A_4$} \\
$\mathfrak{A}_6$ & \multicolumn{2}{c}{$F_{1}=\{1\},F_{2}=\{a,b,d,1\},F_{3}=\{c,d,1\},F_{4}=\{d,1\},F_{5}=A_6$} \\
$\mathfrak{A}_8$ & \multicolumn{2}{c}{$F_{1}=\{1\},F_{2}=\{a,c,d,e,f,1\},F_{3}=\{c,e,1\},F_{4}=\{f,1\},F_{5}=A_8$} \\ \hline
\end{tabular}
\caption{The sets of filters of $\mathfrak{A}_4$, $\mathfrak{A}_6$, and $\mathfrak{A}_8$}
\label{tafiex}
\end{table}
\end{example}

Let $\mathfrak{A}$ and $\mathfrak{B}$ be residuated lattices. A map $f:A\longrightarrow B$ is called a  \textit{morphism of residuated lattices}, in symbols $f:\mathfrak{A}\longrightarrow \mathfrak{B}$, if it preserves the fundamental operations. If $f:\mathfrak{A}\longrightarrow \mathfrak{B}$ is a morphism of residuated lattices, we put $coker(f)=f^{\leftarrow}(1)$. One can see that $f$ is a injective if and only if $coker(f)=\{1\}$.
\begin{example}\label{mora6a4}
  Consider the residuated lattice $\mathfrak{A}_4$ from Example \ref{exa4}, the residuated lattice $\mathfrak{A}_6$ from Example \ref{exa6}. One can see that the map $f:A_{6}\longrightarrow A_{4}$, given by $f(0)=0$, $f(a)=f(b)=a$, $f(c)=b$ and $f(d)=f(1)=1$ is a morphism of residuated lattices.
\end{example}

A proper filter in a residuated lattice $\mathfrak{A}$ is termed \textit{maximal} if it is a maximal element within the set of all proper filters of $\mathfrak{A}$. The collection of maximal filters of $\mathfrak{A}$ is denoted by $Max(\mathfrak{A})$. Additionally, a proper filter $\mathfrak{p}$ in $\mathfrak{A}$ is \textit{prime} if $x\vee y\in \mathfrak{p}$ implies $x\in \mathfrak{p}$ or $y\in \mathfrak{p}$, for any $x, y \in A$. The set of prime filters of $\mathfrak{A}$ is denoted by $Spec(\mathfrak{A})$. Given the distributive lattice property of $\mathscr{F}(\mathfrak{A})$, it follows that $Max(\mathfrak{A})\subseteq Spec(\mathfrak{A})$. Zorn's lemma establishes that any proper filter is contained within a maximal filter and, consequently, within a prime filter. Minimal prime filters, identified as minimal elements within the set of prime filters, form the set $Min(\mathfrak{A})$. For a detailed exploration of prime filters in a residuated lattice, refer to \cite{rasouli2019going}.
\begin{example}\label{maxminex}
Consider the residuated lattice $\mathfrak{A}_4$ from Example \ref{exa4}, the residuated lattice $\mathfrak{A}_6$ from Example \ref{exa6}, and the residuated lattice $\mathfrak{A}_8$ from Example \ref{exa8}. The sets of their maximal, prime, and minimal prime filters are presented in Table \ref{prfiltab}.
\begin{table}[h]
\centering
\begin{tabular}{cccc}
\hline
                 & \multicolumn{3}{c}{Prime filters}      \\ \hline
                 & Maximal filters &          & Minimal prime filters      \\
$\mathfrak{A}_4$ &  $F_{2},F_{3}$                 &  &$F_{2},F_{3}$\\
$\mathfrak{A}_6$ &  $F_2,F_3$                 &  & $F_{1}$\\
$\mathfrak{A}_8$ &  $F_{2}$                          &         & $F_{3},F_{4}$ \\ \hline
\end{tabular}
\caption{The sets of maximal, prime, and minimal prime filters of $\mathfrak{A}_4$, $\mathfrak{A}_6$, and $\mathfrak{A}_8$}
\label{prfiltab}
\end{table}
\end{example}

Let $\mathfrak{A}$ be a residuated lattice, and $\Pi$ a collection of prime filters of $\mathfrak{A}$. For a subset $\pi$ of $\Pi$, we define $k(\pi)=\bigcap \pi$. Additionally, for a subset $X$ of $A$, we set $h_{\Pi}(X)=\{P\in \Pi \mid X\subseteq P\}$ and $d_{\Pi}(X)=\Pi\setminus h_{\Pi}(X)$.

The collection $\Pi$ can be topologized by considering the collection $\{h_{\Pi}(x)\mid x\in A\}$ as a closed (or open) basis. This topological structure is referred to as \textit{the (dual) hull-kernel topology} on $\Pi$ and is denoted by $\Pi_{h(d)}$. Moreover, the topology generated by $\tau_{h}\cup \tau_{d}$ on $\text{Spec}(\mathfrak{A})$ is known as \textit{the patch topology} and is denoted by $\tau_{p}$.

For a detailed discussion of the (dual) hull-kernel and patch topologies on a residuated lattice, we refer to \cite{rasouli2021hull}.

To clarify, we provide some examples of topological spaces associated with residuated lattices.
\begin{example}\label{spectrumofex}
Consider the residuated lattice $\mathfrak{A}_4$ from Example \ref{exa4}, the residuated lattice $\mathfrak{A}_6$ from Example \ref{exa6}, and the residuated lattice $\mathfrak{A}_8$ from Example \ref{exa8}. Table \ref{prfiltab} presents some topological spaces related to these structures.
\begin{table}[h]
\centering
\begin{tabular}{ccl}
\hline
                           & \multicolumn{2}{c}{Topological spaces}      \\ \hline
$Spec_{h}(\mathfrak{A}_4)$ & \multicolumn{2}{c} {$\{\emptyset,\{F_2\},\{F_3\},\{F_2,F_3\}\}$}\\
$Max_{d}(\mathfrak{A}_6)$ & \multicolumn{2}{c} {$\{\emptyset,\{F_{2}\},\{F_{3}\},\{F_2,F_3\}\}$}\\
$Min_{p}(\mathfrak{A}_8)$ & \multicolumn{2}{c} {$\{\emptyset,\{F_{3}\},\{F_{4}\},\{F_{3},F_{4}\},\{F_2,F_3,F_{4}\}\}$} \\ \hline
\end{tabular}
\caption{Some topological spaces associated with $\mathfrak{A}_4$, $\mathfrak{A}_6$, and $\mathfrak{A}_8$}
\label{prfiltab}
\end{table}
\end{example}

\subsection{\'{E}tal\'{e} spaces of sets}

In this subsection, we present key definitions of \'{e}tal\'{e} spaces of sets relevant to our analysis. For completeness, several properties are also established. Additional background can be found in \cite{Bredon1997, Tennison1976}.


A map $f\colon \mathscr{T} \to \mathscr{B}$ is a \emph{local homeomorphism} if, for each point in $\mathscr{T}$, there exists a neighbourhood $U$ such that $f|_U$ is a homeomorphism onto an open subset of $\mathscr{B}$; see \cite[\S 4.4]{engelking1989general}. 

Similarly, a map $f\colon \mathscr{T} \to Y$, where $\mathscr{T}$ is a topological space, is \emph{locally injective} if, for each point in $\mathscr{T}$, there exists a neighbourhood $U$ such that $f|_U$ is injective.


\begin{lemma}\label{lochomeo-properties}
Every local homeomorphism is continuous, open, and locally injective.
\end{lemma}

\begin{proof}
Let $f\colon \mathscr{T} \to \mathscr{B}$ be a local homeomorphism.

To show continuity, let $V \subseteq \mathscr{B}$ be open and fix $t \in f^{-1}(V)$. Since $f$ is a local homeomorphism, there exists a neighbourhood $U$ of $t$ such that $f|_U$ is a homeomorphism onto an open subset of $\mathscr{B}$. Then, $(f|_U)^{-1}(V \cap f(U))$ is open in $U$ and contains $t$, hence $f^{-1}(V)$ is open.

To show openness, let $W \subseteq \mathscr{T}$ be open and fix $b \in f(W)$. Choose $t \in f^{-1}(b) \cap W$. Again, there exists a neighbourhood $U$ of $t$ such that $f|_U$ is a homeomorphism. Then, $f(U \cap W)$ is open in $f(W)$ and contains $b$, hence $f(W)$ is open.

Local injectivity follows directly from the definition of a local homeomorphism, as each local restriction $f|_U$ is a homeomorphism and thus injective.
\end{proof}


\begin{theorem}\label{lhcol}
A map is a local homeomorphism if and only if it is continuous, open, and locally injective. In particular, a bijective local homeomorphism is a homeomorphism.
\end{theorem}

\begin{proof}
The forward direction follows from Lemma~\ref{lochomeo-properties}. Conversely, suppose $f\colon \mathscr{T} \to \mathscr{B}$ is continuous, open, and locally injective. For each $t \in \mathscr{T}$, let $U$ be a neighbourhood such that $f|_U$ is injective. Since $f$ is open and continuous, the image $f(U)$ is open, and $f|_U$ defines a homeomorphism onto $f(U)$. Hence, $f$ is a local homeomorphism.
\end{proof}


\begin{proposition}\label{sheafbasepro}
Let $f\colon \mathscr{T} \to \mathscr{B}$ be a local homeomorphism. Then the family
\[
\mathcal{B} = \{ V \subseteq \mathscr{T} \mid V~\text{is open and}~f|_V~\text{is a homeomorphism} \}
\]
forms a basis for the topology on $\mathscr{T}$.
\end{proposition}

\begin{proof}
Let $U \subseteq \mathscr{T}$ be open and $x \in U$. Since $f$ is a local homeomorphism, there exists an open neighbourhood $W$ of $x$ such that $f|_W$ is a homeomorphism. Then $V := U \cap W$ is open, contains $x$, and satisfies $f|_V$ is a homeomorphism, as it is the restriction of $f|_W$. Hence, $x \in V \subseteq U$ with $V \in \mathcal{B}$, proving that $\mathcal{B}$ is a basis.
\end{proof}


\begin{lemma}\label{lochomeoprop}
\begin{enumerate}
\item [$(1)$ \namedlabel{lochomeoprop1}{$(1)$}] The restriction of a local homeomorphism to any open subspace remains a local homeomorphism;
\item [$(2)$ \namedlabel{lochomeoprop2}{$(2)$}] Any finite composition of local homeomorphisms yields a local homeomorphism.
\end{enumerate}
\end{lemma}

\begin{proof}
The statements follow directly from the definition of local homeomorphisms.
\end{proof}


Let $\mathcal{C}$ be a category and $c$ an object of $\mathcal{C}$. Recall that the \textit{slice category} (a special case of a comma category; see \cite[\S 6, Chapter II]{lane1971categories}) of $\mathcal{C}$ over $c$, denoted by $\mathcal{C}/c$, is the category whose objects are all pairs $(a, f)$, where $a$ is an object of $\mathcal{C}$ and $f$ is a morphism from $a$ to $c$. The morphisms from $(a, f)$ to $(b, g)$ are the morphisms $h$ in $\mathcal{C}$ that make the following diagram commute:
\[
\begin{tikzcd}
a \arrow[rr, "h"] \arrow[rd, "f"'] & & b \arrow[ld, "g"] \\
& c &
\end{tikzcd}
\]
If $\mathcal{C}/c$ is a slice category, then $c$ is referred to as the \textit{base object} of $\mathcal{C}/c$.

In the sequel, \(\textbf{Top}\) denotes the category of topological spaces with continuous maps, while \(\textbf{Top}_{\ell}\) denotes the category of topological spaces with local homeomorphisms as morphisms.

For a given topological space \(\mathscr{B}\), the elements of the slice category \(\mathbf{Etale}(\mathscr{B}) = \textbf{Top}_{\ell}/\mathscr{B}\) are called \textit{\'{e}tal\'{e} spaces of sets over \(\mathscr{B}\)}. The morphisms between them are referred to as \textit{morphisms of \'{e}tal\'{e} spaces of sets over \(\mathscr{B}\)}. If \((\mathscr{T}, \pi)\) is such an object, then \(\mathscr{T}\) is called the \textit{total space}, and \(\pi\) is called the \textit{\'{e}tal\'{e} projection}. For a given point \(b \in \mathscr{B}\), the fiber \(\pi^{-1}(b)\) is referred to as the \textit{stalk of \(\mathscr{T}\) at \(b\)} and is denoted by \(\mathscr{T}_{b}\).

In the following discussion, an \'{e}tal\'{e} space \((\mathscr{T}, \pi)\) of sets over \(\mathscr{B}\) will be denoted by \(\mathscr{T} \overset{\pi}{\downarrow} \mathscr{B}\) and will simply be referred to as an \'{e}tal\'{e} space. Additionally, a morphism of \'{e}tal\'{e} spaces of sets \(h\colon (\mathscr{T}, \pi) \to (\mathscr{S}, \phi)\) over \(\mathscr{B}\) will be denoted by
\[
h\colon \mathscr{T} \overset{\pi}{\downarrow} \mathscr{B} \to \mathscr{S} \overset{\phi}{\downarrow} \mathscr{B}
\]
and referred to as a morphism of \'{e}tal\'{e} spaces. When the base space and the \'{e}tal\'{e} projection are clear from context, they will typically be omitted for brevity.


The next proposition is a restatement of Lemma 3.5 from \cite[\S 2]{Tennison1976}, with a more detailed proof.

\begin{proposition}\label{conopnloc}
Let $\mathscr{T} \overset{\pi}{\downarrow} \mathscr{B}$ and $\mathscr{S} \overset{\phi}{\downarrow} \mathscr{B}$ be two \'{e}tal\'{e} spaces. If $h\colon \mathscr{T} \to \mathscr{S}$ is a map such that $\phi \circ h = \pi$, the following are equivalent:
\begin{enumerate}
\item [$(1)$ \namedlabel{conopnloc1}{$(1)$}] $h$ is a local homeomorphism;
\item [$(2)$ \namedlabel{conopnloc2}{$(2)$}] $h$ is open;
\item [$(3)$ \namedlabel{conopnloc3}{$(3)$}] $h$ is continuous.
\end{enumerate}
\end{proposition}

\begin{proof}
\item[\ref{conopnloc1}$\Rightarrow$\ref{conopnloc2}:] This follows immediately from Lemma~\ref{lochomeo-properties}, as every local homeomorphism is open.

\item[\ref{conopnloc2}$\Rightarrow$\ref{conopnloc3}:]  
By Proposition~\ref{sheafbasepro}, let $V$ be a basic open subset of $\mathscr{S}$ such that $\phi|_V$ is a homeomorphism. Let $x \in h^{-1}(V)$, so that $\pi(x) \in \phi(V)$. Since $\pi$ is a local homeomorphism, there exists a neighbourhood $U$ of $x$ such that $\pi|_U$ is a homeomorphism and $\pi(U) \subseteq \phi(V)$.

Now define $W := h(U) \cap \phi|_V^{-1}(\pi(U))$. Consider any $t \in \pi|_U^{-1}(\phi|_V(W))$. Then there exist $w \in W$ and $u \in U$ such that
\[
\pi(u) = \phi(h(u)) = \phi(w) = \pi(t).
\]
Since $\pi|_U$ is injective, it follows that $u = t$, and hence $t \in h^{-1}(V)$. Thus, $h^{-1}(V)$ is open, and $h$ is continuous.

\item[\ref{conopnloc3}$\Rightarrow$\ref{conopnloc1}:]  
Let $t \in \mathscr{T}$. There exist neighbourhoods $U$ of $t$ and $V$ of $h(t)$ such that $\pi|_U$ and $\phi|_V$ are homeomorphisms. Since $h$ is continuous, the set $W := U \cap h^{-1}(V)$ is an open neighbourhood of $t$.

Suppose $t_1, t_2 \in W$ with $h(t_1) = h(t_2)$. Then
\[
\pi(t_1) = \phi(h(t_1)) = \phi(h(t_2)) = \pi(t_2).
\]
As $\pi|_U$ is injective, it follows that $t_1 = t_2$. Hence, $h|_W$ is injective.

Now let $O \subseteq W$ be open. Then
\[
h(O) = (\phi|_V)^{-1}(\phi(h(O))) = (\phi|_V)^{-1}(\pi(O)),
\]
which is open since $\pi$ and $\phi$ are local homeomorphisms. Thus, $h|_W$ is a homeomorphism onto an open subset, and $h$ is a local homeomorphism.
\end{proof}


\begin{definition}
Let $\mathscr{T} \overset{\pi}{\downarrow} \mathscr{B}$ be an \'{e}tal\'{e} space. A continuous map $\sigma\colon X \to \mathscr{T}$, where $X$ is a subspace of $\mathscr{B}$, is called a \textit{section of $\mathscr{T}$ over $X$} if $\pi \circ \sigma = \mathrm{id}_X$. The set of all such sections is denoted by $\Gamma(X, \mathscr{T})$. Elements of $\Gamma(\mathscr{B}, \mathscr{T})$ are called \textit{global sections} of $\mathscr{T}$.
\end{definition}

\begin{remark}
The nonemptiness of the set of global sections of an \'{e}tal\'{e} space implies that its \'{e}tal\'{e} projection is surjective.
\end{remark}


\begin{lemma}\label{shepropo}
Let $\mathscr{T} \overset{\pi}{\downarrow} \mathscr{B}$ be an \'{e}tal\'{e} space, and let $t \in \mathscr{T}$. Then there exists a neighbourhood $U$ of $\pi(t)$ and a section $\sigma$ of $\mathscr{T}$ over $U$ such that $\sigma(\pi(t)) = t$.
\end{lemma}

\begin{proof}
There exists a neighbourhood $V$ of $t$ such that $\pi|_V$ is a homeomorphism. Set $U := \pi(V)$ and define $\sigma\colon U \to \mathscr{T}$ by $\sigma(x) := (\pi|_V)^{-1}(x)$ for all $x \in U$.
\end{proof}


\begin{lemma}\label{1shepropo}
Let $\mathscr{T} \overset{\pi}{\downarrow} \mathscr{B}$ be an \'{e}tal\'{e} space, let $\sigma$ be a section of $\mathscr{T}$ over an open subset of $\mathscr{B}$, and let $x$ be a point in the domain of $\sigma$. Suppose $V$ is a neighbourhood of $\sigma(x)$ such that $\pi|_V$ is a homeomorphism. Then there exists a neighbourhood $U$ of $x$ such that $\sigma(U) \subseteq V$ and 
\[
\sigma|_U = (\pi|_V)^{-1}|_U.
\]
\end{lemma}

\begin{proof}
By assumption, the set $U := \sigma^{-1}(V)$ is a neighbourhood of $x$ in $\mathscr{B}$. Clearly, $\sigma(U) \subseteq V$. Moreover, for any $x \in U$,
\[
(\pi|_V)^{-1}|_U(x) = (\pi|_V)^{-1}(\pi(\sigma(x))) = \sigma(x),
\]
since $\pi \circ \sigma = \mathrm{id}$ on the domain of $\sigma$, and $(\pi|_V)^{-1}$ is the inverse of $\pi|_V$. Therefore, $\sigma|_U = (\pi|_V)^{-1}|_U$.
\end{proof}


If $X, Y \subseteq \mathscr{B}$ and $\sigma \in \Gamma(X, \mathscr{T})$, $\tau \in \Gamma(Y, \mathscr{T})$, the \textit{equalizer} of $\sigma|_{X \cap Y}$ and $\tau|_{X \cap Y}$ in the category \textbf{Top} is denoted, for simplicity, by $\textsc{Eq}(\sigma, \tau)$. It is defined as the set
\[
\textsc{Eq}(\sigma, \tau) := \{ b \in X \cap Y \mid \sigma(b) = \tau(b) \},
\]
equipped with the subspace topology. By a slight abuse of terminology, we refer to this set as the equalizer of $\sigma$ and $\tau$.


\begin{proposition}\label{sheprop}
Let $\mathscr{T} \overset{\pi}{\downarrow} \mathscr{B}$ be an \'{e}tal\'{e} space. The following statements hold:
\begin{enumerate}
\item [$(1)$ \namedlabel{sheprop1}{$(1)$}] If $U$ and $V$ are open subsets of $\mathscr{B}$, $\sigma \in \Gamma(U, \mathscr{T})$ and $\tau \in \Gamma(V, \mathscr{T})$, then $\textsc{Eq}(\sigma, \tau)$ is open;
\item [$(2)$ \namedlabel{sheprop2}{$(2)$}] If $\mathscr{T}$ is Hausdorff, $U$ and $V$ are open subsets of $\mathscr{B}$, $\sigma \in \Gamma(U, \mathscr{T})$ and $\tau \in \Gamma(V, \mathscr{T})$, then $\textsc{Eq}(\sigma, \tau)$ is clopen;
\item [$(3)$ \namedlabel{sheprop3}{$(3)$}] Each section of $\mathscr{T}$, defined over an open subset of $\mathscr{B}$, is an open map.
\end{enumerate}
\end{proposition}

\begin{proof}
\item[\ref{sheprop1}:] Let $b \in \textsc{Eq}(\sigma, \tau)$. Choose a neighbourhood $W$ of $\sigma(b)$ such that $\pi|_W$ is a homeomorphism. By Lemma~\ref{1shepropo}, there exist neighbourhoods $U_0$ and $V_0$ of $b$ such that $\sigma(U_0), \tau(V_0) \subseteq W$, and $\sigma|_{U_0} = (\pi|_W)^{-1}|_{U_0}$, $\tau|_{V_0} = (\pi|_W)^{-1}|_{V_0}$. Therefore, $\sigma|_{U_0 \cap V_0} = \tau|_{U_0 \cap V_0}$, and hence $U_0 \cap V_0$ is a neighbourhood of $b$ contained in $\textsc{Eq}(\sigma, \tau)$.

\item[\ref{sheprop2}:] This follows from \ref{sheprop1} and the fact that continuous sections into a Hausdorff space yield closed equalizers; see \cite[Theorem 13.13]{willard2012general}.

\item[\ref{sheprop3}:] Let $U \subseteq \mathscr{B}$ be open, and let $\sigma \in \Gamma(U, \mathscr{T})$. For any $t \in \sigma(U)$, write $t = \sigma(b)$ for some $b \in U$. Choose a neighbourhood $W$ of $t$ such that $\pi|_W$ is a homeomorphism. By Lemma~\ref{1shepropo}, there exists a neighbourhood $V$ of $b$ such that $\sigma(V) \subseteq W$ and $\sigma|_V = (\pi|_W)^{-1}|_V$. Hence,
\[
\sigma(V) = (\pi|_W)^{-1}(V),
\]
which is open in $\mathscr{T}$ since $V$ is open in $\mathscr{B}$ and $\pi|_W$ is a homeomorphism. Thus, every point $t \in \sigma(U)$ has a neighbourhood contained in $\sigma(U)$, so $\sigma(U)$ is open.
\end{proof}


\begin{proposition}\label{1sheafbasepro}
Let $\mathscr{T} \overset{\pi}{\downarrow} \mathscr{B}$ be an \'{e}tal\'{e} space. The family
\[
\{ \sigma(U) \mid U~\text{is open in } \mathscr{B},\ \sigma \in \Gamma(U, \mathscr{T}) \}
\]
is a basis for the topology of $\mathscr{T}$.
\end{proposition}

\begin{proof}
Let $V$ be an open subset of $\mathscr{T}$ such that $\pi|_V$ is a homeomorphism. Set $U := \pi(V)$ and define $\sigma \colon U \to \mathscr{T}$ by
\[
\sigma(x) := (\pi|_V)^{-1}(x), \quad \text{for all } x \in U.
\]
Clearly, $\sigma \in \Gamma(U, \mathscr{T})$ and $V = \sigma(U)$. Hence,
\[
\{ V \subseteq \mathscr{T} \mid V~\text{open},\ \pi|_V~\text{is a homeomorphism} \} \subseteq \{ \sigma(U) \mid U~\text{open in } \mathscr{B},\ \sigma \in \Gamma(U, \mathscr{T}) \}.
\]
By Proposition~\ref{sheafbasepro}, the left-hand family forms a basis for the topology of $\mathscr{T}$, so the right-hand family does as well.
\end{proof}


Let $X$ be a set, $\{Y_i\}_{i \in I}$ a collection of topological spaces, and $\{f_i \colon Y_i \to X\}_{i \in I}$ a collection of maps. The \textit{final topology} on $X$ induced by $\{f_i\}_{i \in I}$ is the finest topology on $X$ for which each map $f_i \colon Y_i \to X$ is continuous. It follows that a subset $U \subseteq X$ is open in the final topology if and only if $f_i^{-1}(U)$ is open in $Y_i$ for every $i \in I$.

\begin{proposition}\label{finalshetopo}
Let $\mathscr{T} \overset{\pi}{\downarrow} \mathscr{B}$ be an \'{e}tal\'{e} space. The following hold:
\begin{enumerate}
\item [$(1)$ \namedlabel{finalshetopo1}{$(1)$}] The topology of $\mathscr{T}$ coincides with the final topology induced by the set of sections of $\mathscr{T}$;
\item [$(2)$ \namedlabel{finalshetopo2}{$(2)$}] The topology induced on each stalk of $\mathscr{T}$ is the discrete topology.
\end{enumerate}
\end{proposition}

\begin{proof}
\item[\ref{finalshetopo1}:] Let $V \subseteq \mathscr{T}$. Suppose $V$ is open in the topology of $\mathscr{T}$. Then for every section $\sigma$ of $\mathscr{T}$, the set $\sigma^{-1}(V)$ is open in the domain of $\sigma$, so $V$ is open in the final topology.

Conversely, suppose $V$ is open in the final topology, and let $x \in V$. By Lemma~\ref{shepropo}, there exists a section $\sigma$ such that $\sigma(\pi(x)) = x$. By Proposition~\ref{sheprop}\ref{sheprop3}, the image $\sigma(\sigma^{-1}(V))$ is an open neighbourhood of $x$ contained in $V$. Thus, $V$ is open in the topology of $\mathscr{T}$.

\item[\ref{finalshetopo2}:] Let $x \in \mathscr{B}$. By Lemma~\ref{shepropo}, for any $y \in \mathscr{T}_x$, there exists a section $\sigma$ such that $\sigma(\pi(y)) = y$. It follows that
\[
\sigma(U) \cap \mathscr{T}_x = \{y\}
\]
for some neighbourhood $U$ of $x$. Hence, $\{y\}$ is open in $\mathscr{T}_x$, and the induced topology on $\mathscr{T}_x$ is discrete.
\end{proof}


\section{\'{E}tal\'{e} Spaces of Residuated Lattices}\label{sec3}

In this section, we introduce the concept of \'{e}tal\'{e} spaces of residuated lattices and the morphisms between them. We show that for a given space $\mathscr{B}$, the collection of \'{e}tal\'{e} spaces of residuated lattices over $\mathscr{B}$, together with the morphisms between them, forms a subcategory of $\mathbf{Etale}(\mathscr{B})$.

\begin{definition}\label{rlshprop}
An \'{e}tal\'{e} space $\mathscr{T} \overset{\pi}{\downarrow} \mathscr{B}$ is called an \emph{\'{e}tal\'{e} space of residuated lattices} if the following conditions hold:
\begin{enumerate}
    \item [$(1)$ \namedlabel{rlshprop1}{$(1)$}] For each point $b \in \mathscr{B}$, the stalk $\mathscr{T}_{b}$ is a residuated lattice. The operations of the residuated lattice $\mathscr{T}_{b}$ are indexed by $b$.

    \item [$(2)$ \namedlabel{rlshprop2}{$(2)$}] For any subset $X \subseteq \mathscr{B}$, the set of sections $\Gamma(X, \mathscr{T})$ forms a residuated lattice, with operations defined pointwise as follows:

    \[
    \begin{array}{lcll}
        \diamond_{X}\colon & \Gamma(X, \mathscr{T}) \times \Gamma(X, \mathscr{T}) & \to & \Gamma(X, \mathscr{T}) \\
        & (\sigma, \rho) & \mapsto & \left( x \mapsto \sigma(x) \diamond_{x} \rho(x) \right)
    \end{array}
    \]
    where $\diamond$ is a binary fundamental operation of residuated lattices; and

    \[
    \begin{array}{llll}
        \diamond_{X}\colon & X & \to & \mathscr{T} \\
        & x & \mapsto & \diamond_{x}
    \end{array}
    \]
    where $\diamond$ is a nullary fundamental operation of residuated lattices.
\end{enumerate}
\end{definition}


Let $\mathscr{T} \overset{\pi}{\downarrow} \mathscr{B}$ be an \'{e}tal\'{e} space of residuated lattices. Define
\[
\mathscr{T}^{(2)} = \{(a_1, a_2) \in \mathscr{T}^2 \mid \pi(a_1) = \pi(a_2)\},
\]
and equip it with the subspace topology inherited from the product topology on $\mathscr{T}^2$. If $\diamond$ is a binary fundamental operation of residuated lattices, and $\diamond_b$ denotes the corresponding operation on the stalk $\mathscr{T}_b$ for each $b \in \mathscr{B}$, then we define a map
\[
\overset{\mathscr{T}}{\diamond} \colon \mathscr{T}^{(2)} \to \mathscr{T}, \quad \text{by} \quad \overset{\mathscr{T}}{\diamond}(a_1, a_2) := \diamond_b(a_1, a_2),
\]
where $\pi(a_1) = \pi(a_2) = b$. Likewise, if $\diamond$ is a nullary function symbol and $\diamond_b$ is the corresponding constant in the stalk $\mathscr{T}_b$, we define
\[
\overset{\mathscr{T}}{\diamond} \colon \mathscr{B} \to \mathscr{T}, \quad \text{by} \quad \overset{\mathscr{T}}{\diamond}(b) := \diamond_b.
\]
In what follows, we will omit the superscript $\mathscr{T}$ when the context is clear.


The next result, inspired by a similar statement for universal algebras \cite[Lemma~1.2]{davey1973sheaf}, provides necessary and sufficient conditions for an \'{e}tal\'{e} space to be an \'{e}tal\'{e} space of residuated lattices.

\begin{proposition}\label{sheafalgeblem}
Let $\mathscr{T} \overset{\pi}{\downarrow} \mathscr{B}$ be an \'{e}tal\'{e} space such that each stalk $\mathscr{T}_b$ is a residuated lattice. The following conditions are equivalent:
\begin{enumerate}
  \item [$(1)$ \namedlabel{sheafalgeblem1}{$(1)$}] $\mathscr{T}$ is an \'{e}tal\'{e} space of residuated lattices over the topological space $\mathscr{B}$;
  \item [$(2)$ \namedlabel{sheafalgeblem2}{$(2)$}] For every open subset $U \subseteq \mathscr{B}$, the set $\Gamma(U, \mathscr{T})$ of sections over $U$ forms a residuated lattice under pointwise operations;
  \item [$(3)$ \namedlabel{sheafalgeblem3}{$(3)$}] The operations $\vee$, $\wedge$, $\odot$, $\to$, $0$, and $1$ are continuous.
\end{enumerate}
\end{proposition}

\begin{proof}
\item[\ref{sheafalgeblem1}$\Rightarrow$\ref{sheafalgeblem2}:] This follows directly from Definition~\ref{rlshprop}.

\item[\ref{sheafalgeblem2}$\Rightarrow$\ref{sheafalgeblem3}:]
Let $U \subseteq \mathscr{B}$ be open and $\sigma \in \Gamma(U, \mathscr{T})$. By Proposition~\ref{1sheafbasepro}, the image $\sigma(U)$ is a basic open set in $\mathscr{T}$. Let $(x, y) \in \vee^{-1}(\sigma(U))$. Then
\[
x \vee_b y = \sigma(b), \quad \text{where } \pi(x) = \pi(y) = b \in U.
\]
By Lemma~\ref{shepropo}, there exist neighbourhoods $W$ of $b$ and sections $\sigma_x, \sigma_y \in \Gamma(W, \mathscr{T})$ such that $\sigma_x(b) = x$ and $\sigma_y(b) = y$. Then $(\sigma_x \vee \sigma_y)(b) = \sigma(b)$, and by Proposition~\ref{sheprop}\ref{sheprop1}, the set
\[
V = \mathrm{Eq}(\sigma, \sigma_x \vee \sigma_y)
\]
is a neighbourhood of $b$. Therefore,
\[
(x, y) \in \sigma_x(V) \times \sigma_y(V) \cap \mathscr{T}^{(2)} \subseteq \vee^{-1}(\sigma(U)),
\]
showing that $\vee$ is continuous. The continuity of $\wedge$, $\odot$, and $\to$ follows by analogous arguments. The continuity of $0$ and $1$ is clear, since they are defined pointwise via $b \mapsto 0_b$ and $b \mapsto 1_b$.

\item[\ref{sheafalgeblem3}$\Rightarrow$\ref{sheafalgeblem1}:]
Let $X \subseteq \mathscr{B}$ and $\sigma, \tau \in \Gamma(X, \mathscr{T})$. For any binary fundamental operation $\diamond$, the pointwise operation $\diamond_X(\sigma, \tau)$ is given by
\[
\diamond_X(\sigma, \tau) = \diamond \circ (\sigma, \tau),
\]
which is continuous, since $\sigma$, $\tau$, and $\diamond$ are continuous. Hence, $\diamond_X(\sigma, \tau) \in \Gamma(X, \mathscr{T})$.

For each nullary operation symbol $\diamond$, define the constant section $\diamond_X \colon X \to \mathscr{T}$ by
\[
\diamond_X = \diamond|_X.
\]
Then $(\Gamma(X, \mathscr{T}); \vee_X, \wedge_X, \odot_X, \to_X, 0_X, 1_X)$ forms a residuated lattice under pointwise operations. Therefore, $\mathscr{T}$ satisfies Definition~\ref{rlshprop}.
\end{proof}


\begin{definition}\label{resetalmor}
Let $\mathscr{T} \overset{\pi}{\downarrow} \mathscr{B}$ and $\mathscr{S} \overset{\phi}{\downarrow} \mathscr{B}$ be two \'{e}tal\'{e} spaces of residuated lattices. A morphism of \'{e}tal\'{e} spaces $h\colon \mathscr{T} \to \mathscr{S}$ is called a \emph{morphism of \'{e}tal\'{e} spaces of residuated lattices} if, for every point $b \in \mathscr{B}$, the restriction $h|_{\mathscr{T}_b}$—denoted $h_b$—is a morphism of residuated lattices.
\end{definition}


\begin{remark}
Let $\mathscr{T} \overset{\pi}{\downarrow} \mathscr{B}$ and $\mathscr{S} \overset{\phi}{\downarrow} \mathscr{B}$ be two \'{e}tal\'{e} spaces of residuated lattices. To define a morphism of \'{e}tal\'{e} spaces of residuated lattices from $\mathscr{T}$ to $\mathscr{S}$, it suffices to specify a family $\{h_b \colon \mathscr{T}_b \to \mathscr{S}_b\}_{b \in \mathscr{B}}$ of residuated lattice morphisms. The global morphism $h$ is then defined by
\[
h(t) := h_{\pi(t)}(t), \quad \text{for all } t \in \mathscr{T}.
\]
Moreover, $h$ is injective (respectively, surjective) if and only if each $h_b$ is injective (respectively, surjective) for all $b \in \mathscr{B}$.
\end{remark}


\begin{theorem}\label{etressubcatetsp}
For a given topological space $\mathscr{B}$, the class of \'{e}tal\'{e} spaces of residuated lattices over $\mathscr{B}$, together with the morphisms of \'{e}tal\'{e} spaces of residuated lattices, forms a subcategory of $\mathbf{Etale}(\mathscr{B})$. This subcategory is denoted by $\mathbf{RL\text{-}Etale}(\mathscr{B})$.
\end{theorem}

\begin{proof}
The result follows from the observation that identity morphisms and compositions of morphisms of \'{e}tal\'{e} spaces of residuated lattices over $\mathscr{B}$ remain morphisms of such spaces. Specifically, for each $b \in \mathscr{B}$, the identity map on $\mathscr{T}_b$ is a residuated lattice morphism, and the composition of two such morphisms $h_b \colon \mathscr{T}_b \to \mathscr{S}_b$ and $k_b \colon \mathscr{S}_b \to \mathscr{R}_b$ yields another residuated lattice morphism $k_b \circ h_b \colon \mathscr{T}_b \to \mathscr{R}_b$. Thus, the class of objects and morphisms satisfies the axioms of a subcategory.
\end{proof}


\section{Presheaves and Sheaves of Residuated Lattices}\label{sec4}

This section discusses presheaves and sheaves of residuated lattices. It introduces the concepts of a presheaf and a sheaf of residuated lattices on a topological space and provides several examples of these notions. Finally, a functor from the category of \'{e}tal\'{e} spaces of residuated lattices to the category of sheaves of residuated lattices is constructed.


\begin{definition}
A pair $(\mathscr{F}, \mathscr{B})$ is called a \textit{presheaf of residuated lattices} if $\mathscr{B}$ is a topological space and $\mathscr{F}$ is a contravariant functor from the category $\mathscr{O}(\mathscr{B})$ of open subsets of $\mathscr{B}$ (with inclusions as morphisms) to the category $\mathbf{RL}$ of residuated lattices. That is,
\[
\mathscr{F} \colon \mathscr{O}(\mathscr{B})^{\mathrm{op}} \to \mathbf{RL}.
\]
\end{definition}


Let $(\mathscr{F}, \mathscr{B})$ be a presheaf of residuated lattices. For every open set $U \subseteq \mathscr{B}$, the elements of $\mathscr{F}(U)$ are referred to as the \textit{sections of $\mathscr{F}$ over $U$}. The elements of $\mathscr{F}(\mathscr{B})$ are called the \textit{global sections} of $\mathscr{F}$. 

If $U \subseteq V$, the corresponding morphism of residuated lattices $\mathscr{F}(V) \to \mathscr{F}(U)$ is denoted by $\mathscr{F}_{V,U}$. In analogy with the restriction of functions, we often write $\mathscr{F}_{V,U}(s)$ as $s|_U$ and refer to it as the \textit{restriction of $s$ to $U$}.

In what follows, when we say that $\mathscr{F}$ is a \textit{presheaf of residuated lattices over} $\mathscr{B}$, we mean that the pair $(\mathscr{F}, \mathscr{B})$ is a presheaf of residuated lattices. The space $\mathscr{B}$ is referred to as the \textit{base space} of the presheaf, and when it is clear from context, it will typically be omitted.


\begin{example}\label{onepointpre}
Let $\mathfrak{A}$ be a residuated lattice and let $\mathscr{B}$ be a one-point topological space. Define a presheaf $\mathscr{P}^{\mathfrak{A}}$ by setting:
\[
\mathscr{P}^{\mathfrak{A}}(\mathscr{B}) = \mathfrak{A}, \quad \text{and} \quad \mathscr{P}^{\mathfrak{A}}(\emptyset) = \{1\}.
\]
It is straightforward to verify that $\mathscr{P}^{\mathfrak{A}}$ defines a presheaf of residuated lattices.
\end{example}


\begin{example}\label{skyscraperpre}
Let $\mathfrak{A}$ be a residuated lattice and let $\mathscr{B}$ be a topological space. Fix a point $b \in \mathscr{B}$. Define a presheaf of residuated lattices $\mathscr{S}^{b, \mathfrak{A}}$ by:
\[
\mathscr{S}^{b, \mathfrak{A}}(U) =
\begin{cases}
\mathfrak{A}, & \text{if } b \in U, \\
\{1\}, & \text{if } b \notin U,
\end{cases}
\]
for each open set $U \subseteq \mathscr{B}$. The restriction maps are defined naturally as follows: if $U \subseteq V$ are open subsets of $\mathscr{B}$ and $b \in V$, define
\[
\mathscr{S}^{b, \mathfrak{A}}_{V,U}(x) =
\begin{cases}
x, & \text{if } b \in U, \\
1, & \text{if } b \notin U,
\end{cases}
\quad \text{for } x \in \mathscr{S}^{b, \mathfrak{A}}(V).
\]
This presheaf is called the \textit{skyscraper presheaf of residuated lattices} at $b$ with value $\mathfrak{A}$.
\end{example}


\begin{example}\label{constantpreex}
Let $\mathfrak{A}$ be a residuated lattice and $\mathscr{B}$ a topological space. A \textit{constant presheaf of residuated lattices} $\mathscr{C}^{\mathfrak{A}}$ is defined by setting
\[
\mathscr{C}^{\mathfrak{A}}(U) = \mathfrak{A}
\]
for every open set $U \subseteq \mathscr{B}$. The restriction maps are taken to be $\mathscr{C}^{\mathfrak{A}}_{V,U} = \mathrm{id}_{\mathfrak{A}}$ for all open sets $U \subseteq V$ in $\mathscr{B}$.
\end{example}


\begin{example}\label{fundexampresh}
Let $\mathfrak{A}$ be a residuated lattice, and let $\{F_b\}_{b \in \mathscr{B}}$ be a family of proper filters of $\mathfrak{A}$ indexed by a topological space $\mathscr{B}$. Define a presheaf $\mathscr{F}^{\mathfrak{A}}$ by
\[
\mathscr{F}^{\mathfrak{A}}(U) := \mathfrak{A} / F_U, \quad \text{where } F_U := \bigcap_{b \in U} F_b,
\]
for every open set $U \subseteq \mathscr{B}$. If $U \subseteq V$, then $F_V \subseteq F_U$, and the restriction map
\[
\mathscr{F}^{\mathfrak{A}}_{V,U} \colon \mathscr{F}^{\mathfrak{A}}(V) \to \mathscr{F}^{\mathfrak{A}}(U)
\]
is defined by
\[
\mathscr{F}^{\mathfrak{A}}_{V,U}(a/F_V) := a/F_U,
\]
which is well-defined due to the inclusion $F_V \subseteq F_U$.
\end{example}


\begin{example}\label{prsresexa4}
Consider the residuated lattice $\mathfrak{A}_4$ from Example~\ref{exa4}. By Example~\ref{spectrumofex}, we have
\[
\mathrm{Spec}_h(\mathfrak{A}_4) = \{\emptyset, U = \{F_2\}, V = \{F_3\}, B = \{F_2, F_3\}\}.
\]
Define $\mathscr{F} \colon \mathrm{Spec}_h(\mathfrak{A}_4) \to \mathbf{RL}$ by
\[
\mathscr{F}(\emptyset) = \mathfrak{A}_4 / F_4, \quad \mathscr{F}(U) = \mathfrak{A}_4 / F_2, \quad \mathscr{F}(V) = \mathfrak{A}_4 / F_3, \quad \mathscr{F}(B) = \mathfrak{A}_4 / F_1.
\]
Then $\mathscr{F}$ is a presheaf of residuated lattices.
\end{example}


\begin{example}\label{prsresexa6}
Consider the residuated lattice $\mathfrak{A}_6$ from Example~\ref{exa6}. By Example~\ref{spectrumofex}, we have
\[
\mathrm{Max}_d(\mathfrak{A}_6) = \{\emptyset, U = \{F_2\}, V = \{F_3\}, B = \{F_2, F_3\}\}.
\]
Define $\mathscr{F} \colon \mathrm{Max}_d(\mathfrak{A}_6) \to \mathbf{RL}$ by
\[
\mathscr{F}(\emptyset) = \mathfrak{A}_6 / F_5, \quad \mathscr{F}(U) = \mathfrak{A}_6 / F_2, \quad \mathscr{F}(V) = \mathfrak{A}_6 / F_3, \quad \mathscr{F}(B) = \mathfrak{A}_6 / F_1.
\]
Then $\mathscr{F}$ is a presheaf of residuated lattices.
\end{example}


\begin{example}\label{prsresexa8}
Consider the residuated lattice $\mathfrak{A}_8$ from Example~\ref{exa8}. By Example~\ref{spectrumofex}, we have
\[
\mathrm{Min}_p(\mathfrak{A}_8) = \{\emptyset, U = \{F_3\}, V = \{F_4\}, W = \{F_3, F_4\}, B = \{F_2, F_3, F_4\}\}.
\]
Define $\mathscr{F} \colon \mathrm{Min}_p(\mathfrak{A}_8) \to \mathbf{RL}$ by
\[
\begin{array}{lll}
  \mathscr{F}(\emptyset) = \mathfrak{A}_8 / F_5 & \mathscr{F}(U) = \mathfrak{A}_8 / F_3 & \mathscr{F}(V) = \mathfrak{A}_8 / F_4\\
  \mathscr{F}(W) = \mathfrak{A}_8 / F_1 & \mathscr{F}(B) = \mathfrak{A}_8 / F_1 &
\end{array}
\]
Then $\mathscr{F}$ is a presheaf of residuated lattices.
\end{example}


The following example establishes a connection between topology, algebra, and logic by employing a presheaf to model context-dependent fuzzy truth values. It serves as a foundational case study for more advanced applications, including the logical semantics of spatial and temporal systems, as well as multi-context reasoning frameworks.

\begin{example}\label{sirpshe}
Let \( \mathscr{B} = \{x, y\} \) denote the Sierpi\'{n}ski space endowed with the topology \( \mathscr{O}(\mathscr{B}) = \{\emptyset, \{x\}, \{x, y\}\} \). We define a presheaf \( \mathscr{F} \) of residuated lattices over \( \mathscr{B} \) as follows.

For each open set \( U \subseteq \mathscr{B} \), the set \( \mathscr{F}(U) \) consists of fuzzy truth assignments over the propositional variables \( \{\alpha, \beta\} \), where the structure of \( \mathscr{F}(U) \) depends on the nature of the open set \( U \). For instance, we may interpret \( \alpha \) as the proposition ``It is raining at location \( x \)" and \( \beta \) as the proposition ``It is sunny at location \( y \)".

The set \( \mathscr{F}(\{x, y\}) \) comprises all fuzzy truth assignments for the propositions \( \alpha \) and \( \beta \), that is, all functions \( \nu: \{\alpha, \beta\} \to [0,1] \). These assignments form a pointwise residuated lattice, where the partial order is given by
\[
\nu \leq \mu \quad \text{if and only if} \quad \nu(\chi) \leq \mu(\chi) \quad \text{for all } \chi \in \{\alpha, \beta\},
\]
the monoidal product is defined by
\[
(\nu \odot \mu)(\chi) = \max\{0, \nu(\chi) + \mu(\chi) - 1\},
\]
and the residuum is given by
\[
(\nu \to \mu)(\chi) = \min\{1, 1 - \nu(\chi) + \mu(\chi)\}.
\]

The set \( \mathscr{F}(\{x\}) \) consists of fuzzy assignments restricted to the proposition \( \alpha \), namely, functions \( \nu: \{\alpha\} \to [0,1] \), with \( \beta \) excluded as it is not observable at the point \( x \). The algebraic operations are defined analogously, restricted to the domain \( \{\alpha\} \).

Finally, we define \( \mathscr{F}(\emptyset) = \{\top\} \), a singleton set representing the trivial residuated lattice, where \( \top \) denotes the unique (empty) assignment. This reflects the intuition that, in the absence of any observable data, no meaningful truth assignment can be specified.

The restriction maps describe how truth assignments are localized from larger to smaller open sets. The restriction map \( \mathscr{F}_{\{x,y\}\{x\}} \) restricts a section \( \nu \in \mathscr{F}(\{x,y\}) \), defined on \( \{\alpha, \beta\} \), to its component on \( \{\alpha\} \). In other words, the restriction omits any reference to \( \beta \), which is not observable in the smaller open set \( \{x\} \). The restriction maps \( \mathscr{F}_{\{x, y\}\emptyset} \) and \( \mathscr{F}_{\{x\}\emptyset} \) map any assignment to \( \top \), thereby collapsing all information to the unique element of \( \mathscr{F}(\emptyset) \).

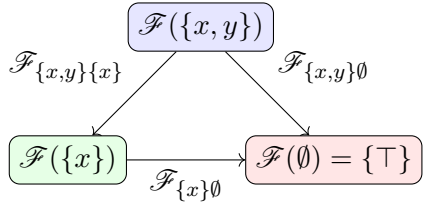
\begin{center}
\begin{tikzpicture}[node distance=2.5cm, auto]
  \node (Uxy) [draw, rounded corners, fill=blue!10] {\( \mathscr{F}(\{x,y\}) \)};
  \node (Ux) [below left of=Uxy, draw, rounded corners, fill=green!10] {\( \mathscr{F}(\{x\}) \)};
  \node (Uempty) [below right of=Uxy, draw, rounded corners, fill=red!10] {\( \mathscr{F}(\emptyset) = \{\top\} \)};

  \draw[->] (Uxy) to node [swap] {\( \mathscr{F}_{\{x,y\}\{x\}} \)} (Ux);
  \draw[->] (Uxy) to node {\( \mathscr{F}_{\{x,y\}\emptyset} \)} (Uempty);
  \draw[->] (Ux) to node [swap] {\( \mathscr{F}_{\{x\}\emptyset} \)} (Uempty);
\end{tikzpicture}
\end{center}

The presheaf axioms are readily verified. For any open set \( U \subseteq \mathscr{B} \), the identity condition \( \mathscr{F}_{UU} = \mathrm{id}_{\mathscr{F}(U)} \) holds, as restricting a function to the same domain leaves it unchanged. The compositionality of restriction maps follows directly from their definitions. For example, given any \( \nu \in \mathscr{F}(\{x, y\}) \), we have
\[
\mathscr{F}_{\{x\}\emptyset}(\mathscr{F}_{\{x, y\}\{x\}}(\nu)) = \mathscr{F}_{\{x\}\emptyset}(\nu|_{\{\alpha\}}) = \top = \mathscr{F}_{\{x, y\}\emptyset}(\nu),
\]
which verifies the required compatibility condition. All other instances of compositionality can be verified analogously.

Hence, the construction \( \mathscr{F} \) defines a valid presheaf of residuated lattices over the Sierpi\'{n}ski space \( \mathscr{B} \).

\end{example}


\begin{theorem}
The category of presheaves of residuated lattices over a topological space $\mathscr{B}$, with natural transformations as morphisms and composition given by vertical composition, forms a category. This category is denoted by $\mathbf{RL\text{-}PreSheaf}(\mathscr{B})$.
\end{theorem}

\begin{proof}
The result is straightforward. The identity natural transformation acts as the identity morphism, and vertical composition of natural transformations is associative. Thus, the category axioms are satisfied.
\end{proof}


The concept of a presheaf is valuable because it enables us to work with objects that are defined only locally rather than globally. Later on, we introduce the concept of a sheaf, which allows us to glue together locally defined data to obtain global objects defined over unions of open sets.

\begin{definition}\label{pressheadef}
Let $\mathscr{F}$ be a presheaf of residuated lattices over a topological space $\mathscr{B}$. Then:
\begin{enumerate}
  \item [$(S)$ \namedlabel{s}{$(S)$}] $\mathscr{F}$ satisfies the \textit{separation property} if, for any open set $O \subseteq \mathscr{B}$, any open cover $\mathscr{U}$ of $O$, and any two sections $s, t \in \mathscr{F}(O)$ such that $s|_U = t|_U$ for all $U \in \mathscr{U}$, it follows that $s = t$.

  \item [$(G)$ \namedlabel{g}{$(G)$}] $\mathscr{F}$ satisfies the \textit{gluing property} if, for any open set $O \subseteq \mathscr{B}$, any open cover $\mathscr{U}$ of $O$, and any family $\{s_U \in \mathscr{F}(U) \mid U \in \mathscr{U}\}$ such that $s_U|_{U \cap V} = s_V|_{U \cap V}$ for all $U, V \in \mathscr{U}$, there exists a section $s \in \mathscr{F}(O)$ such that $s|_U = s_U$ for all $U \in \mathscr{U}$.
\end{enumerate}
\end{definition}


\begin{proposition}
Let $\mathscr{F}$ be a presheaf of residuated lattices over $\mathscr{B}$. Then $\mathscr{F}$ satisfies property \ref{s} if and only if the following holds: for any open set $O \subseteq \mathscr{B}$, any open cover $\mathscr{U}$ of $O$, and any section $s \in \mathscr{F}(O)$ such that $s|_U = 1$ for all $U \in \mathscr{U}$, we have $s = 1$. In other words, a section is globally the unit element if it is locally the unit element.
\end{proposition}

\begin{proof}
Assume that $\mathscr{F}$ satisfies property \ref{s}. Let $O \subseteq \mathscr{B}$ be open, and let $\mathscr{U}$ be an open cover of $O$. Suppose $s \in \mathscr{F}(O)$ satisfies $s|_U = 1$ for every $U \in \mathscr{U}$. Then $s|_U = 1|_U$ for all $U \in \mathscr{U}$, so by separation, we conclude $s = 1$.

Conversely, suppose the stated unit condition holds. Let $s, t \in \mathscr{F}(O)$ be such that $s|_U = t|_U$ for every $U \in \mathscr{U}$. Then for each $U \in \mathscr{U}$, we have
\[
(s \rightarrow t)|_U = s|_U \rightarrow t|_U = 1, \quad \text{and} \quad (t \rightarrow s)|_U = t|_U \rightarrow s|_U = 1.
\]
So both $s \rightarrow t$ and $t \rightarrow s$ restrict to the unit section on each $U \in \mathscr{U}$, and by hypothesis, this implies $s \rightarrow t = 1$ and $t \rightarrow s = 1$, whence $s = t$.
\end{proof}


The following result provides a concise characterization of properties \ref{s} and \ref{g}.

\begin{theorem}\label{preshshesg}
Let $\mathscr{F}$ be a presheaf of residuated lattices over a topological space $\mathscr{B}$. Then $\mathscr{F}$ is a sheaf of residuated lattices if and only if, for every open set $O \subseteq \mathscr{B}$ and every open cover $\mathscr{U}$ of $O$, the following diagram is an equalizer:
\[
\begin{tikzcd}
\mathscr{F}(O) \arrow[rr, "f" description] & & \displaystyle\prod_{U \in \mathscr{U}} \mathscr{F}(U)
\arrow[r, "g" description, bend left=49]
\arrow[r, "h" description, bend right=49] & \displaystyle\prod_{(U, V) \in \mathscr{U}^2} \mathscr{F}(U \cap V)
\end{tikzcd}
\]
where:
\begin{itemize}
  \item $f(s) = (s|_U)_{U \in \mathscr{U}}$;
  \item $g((s_U)_{U \in \mathscr{U}}) = (s_U|_{U \cap V})_{(U, V) \in \mathscr{U}^2}$;
  \item $h((s_U)_{U \in \mathscr{U}}) = (s_V|_{U \cap V})_{(U, V) \in \mathscr{U}^2}$.
\end{itemize}
\end{theorem}

\begin{proof}
This is a direct reformulation of Definition~\ref{pressheadef}, where the equalizer condition captures both the separation and gluing properties.
\end{proof}


\begin{definition}
A presheaf $\mathscr{F}$ of residuated lattices over a topological space $\mathscr{B}$ is called a \textit{sheaf} if it satisfies both the separation property \ref{s} and the gluing property \ref{g}. The collection of sheaves of residuated lattices over $\mathscr{B}$ forms a full subcategory of $\mathbf{RL\text{-}PreSheaf}(\mathscr{B})$, denoted by $\mathbf{RL\text{-}Sheaf}(\mathscr{B})$.
\end{definition}


\begin{example}\label{onepointshe}
The presheaves of residuated lattices introduced in Examples~\ref{onepointpre} and~\ref{skyscraperpre} satisfy the sheaf condition, and thus define sheaves of residuated lattices.
\end{example}


\begin{example}\label{sirppresheaf}
It is straightforward to verify that the presheaf of residuated lattices described in Example~\ref{sirpshe} satisfies the sheaf condition, and thus constitutes a sheaf of residuated lattices.
\end{example}


The following examples illustrates a presheaf of residuated lattices that is not a sheaf.

\begin{example}\label{exsheanot}
Consider the constant presheaf of residuated lattices $\mathscr{C}^{\mathfrak{A}}$ from Example~\ref{constantpreex}. It is straightforward to verify that $\mathscr{C}^{\mathfrak{A}}$ satisfies the separation property \ref{s}. However, in general, it does not satisfy the gluing property \ref{g}.

To see this, let $U_1$ and $U_2$ be disjoint open subsets of $\mathscr{B}$, and let $U := U_1 \cup U_2$. Let $s_1 \in \mathscr{C}^{\mathfrak{A}}(U_1)$ and $s_2 \in \mathscr{C}^{\mathfrak{A}}(U_2)$ with $s_1 \neq s_2$. Since $\mathscr{C}^{\mathfrak{A}}(U) = \mathfrak{A}$ and the restriction maps are identities, there is no $s \in \mathscr{C}^{\mathfrak{A}}(U)$ such that $s|_{U_1} = s_1$ and $s|_{U_2} = s_2$. Hence, the gluing condition fails.
\end{example}


\begin{example}\label{exsirsheanot}
Consider a variation of the presheaf \( \mathscr{F} \) described in Example~\ref{sirpshe}, in which the definition of \( \mathscr{F}(\{x\}) \) is modified. Instead of omitting the proposition \( \beta \), we now define \( \mathscr{F}(\{x\}) \) to consist of all fuzzy truth assignments \( \nu: \{\alpha, \beta\} \to [0,1] \) such that \( \nu(\beta) = 0 \). In this setting, \( \beta \) is not excluded from consideration, but is instead explicitly constrained to have truth value zero on the open set \( \{x\} \). This models a context in which \( \beta \) is observable at location \( x \), but is definitively false.

Accordingly, the restriction map \( \mathscr{F}_{\{x,y\}\{x\}} \) is modified to enforce this constraint: for any \( \nu \in \mathscr{F}(\{x,y\}) \) with \( \nu(\alpha) = a \) and \( \nu(\beta) = b \), we define
\[
\mathscr{F}_{\{x,y\}\{x\}}(\nu) = (\alpha \mapsto a,\; \beta \mapsto 0).
\]
All other restriction maps remain as in the original construction; in particular, restrictions to the empty set are trivial, mapping every assignment to the unique element \( \top \in \mathscr{F}(\emptyset) \).

This modified presheaf fails to satisfy the sheaf condition. Consider the family of sections given by
\[
\nu = (\alpha \mapsto 0.5,\; \beta \mapsto 0) \in \mathscr{F}(\{x\}), \quad \mu = (\alpha \mapsto 0.5,\; \beta \mapsto 0.7) \in \mathscr{F}(\{x,y\}).
\]
These sections are compatible in the sense that
\[
\mathscr{F}_{\{x,y\}\{x\}}(\mu) = (\alpha \mapsto 0.5,\; \beta \mapsto 0) = \nu.
\]
However, there exists no global section \( \sigma \in \mathscr{F}(\{x,y\}) \) that restricts to both \( \nu \) and \( \mu \). Indeed, any such \( \sigma \) would have to equal \( \mu \), but \( \mu \) does not restrict to \( \nu \) under the modified restriction map—since \( \mu(\beta) = 0.7 \), and the restriction forces \( \beta \mapsto 0 \) on \( \{x\} \). Thus, the gluing fails.

The obstruction arises from the imposed constraint that all sections over \( \{x\} \) must assign the value zero to the proposition \( \beta \), thereby preventing gluing with any section over \( \{x,y\} \) in which \( \beta \) assumes a nonzero value.
\end{example}


The rest of this section is devoted to constructing a functor from the category $\textbf{RL-Etale}(\mathscr{B})$ to the category $\textbf{RL-PreSheaf}(\mathscr{B})$.

\begin{proposition}\label{contavafunctshe}
Let $\mathscr{T} \overset{\pi}{\downarrow} \mathscr{B}$ be an \'{e}tal\'{e} space of residuated lattices. Define the map
\[
\textsc{Ps}(\mathscr{T}) \colon \mathscr{O}(\mathscr{B})^{\mathrm{op}} \to \mathbf{RL}
\]
by assigning to each open set \( U \subseteq \mathscr{B} \) the residuated lattice \( \Gamma(U, \mathscr{T}) \), and to each inclusion \( U \subseteq V \) the restriction map
\[
\textsc{Ps}(\mathscr{T})_{V,U} \colon \Gamma(V, \mathscr{T}) \to \Gamma(U, \mathscr{T}), \quad \sigma \mapsto \sigma|_U.
\]
Then \(\textsc{Ps}(\mathscr{T})\) is a contravariant functor.
\end{proposition}

\begin{proof}
This follows directly from the functorial properties of restrictions: identity maps restrict to identity, and composition of restrictions corresponds to restriction over nested open sets.
\end{proof}


\begin{definition}
Let $\mathscr{T} \overset{\pi}{\downarrow} \mathscr{B}$ be an \'{e}tal\'{e} space of residuated lattices. The presheaf \(\textsc{Ps}(\mathscr{T})\) defined above is called the \textit{presheaf of sections} of $\mathscr{T}$.
\end{definition}


\begin{theorem}\label{etaleshea}
Let $\mathscr{T} \overset{\pi}{\downarrow} \mathscr{B}$ be an \'{e}tal\'{e} space of residuated lattices. Then the presheaf of sections \(\textsc{Ps}(\mathscr{T})\) is a sheaf of residuated lattices.
\end{theorem}

\begin{proof}
Let \( O \in \mathscr{O}(\mathscr{B}) \) and let \(\mathscr{U}\) be an open cover of \(O\).

To verify the separation property \ref{s}, suppose \( s, t \in \Gamma(O, \mathscr{T}) \) satisfy \( s|_U = t|_U \) for all \( U \in \mathscr{U} \). For any \( x \in O \), there exists \( U \in \mathscr{U} \) such that \( x \in U \), and thus:
\[
s(x) = s|_U(x) = t|_U(x) = t(x).
\]
Hence, \( s = t \), so the presheaf satisfies the separation property.

To verify the gluing property \ref{g}, suppose we are given a family \(\{ s_U \in \Gamma(U, \mathscr{T}) \mid U \in \mathscr{U} \}\) such that
\[
s_U|_{U \cap V} = s_V|_{U \cap V} \quad \text{for all } U, V \in \mathscr{U}.
\]
Define a function \( s \colon O \to \mathscr{T} \) by setting \( s(x) := s_U(x) \) for any \( U \in \mathscr{U} \) such that \( x \in U \). This is well-defined by the compatibility condition on overlaps.

To show continuity of \( s \), note that for any \( x \in O \), there exists \( U \in \mathscr{U} \) such that \( x \in U \), and \( s \) agrees with the continuous section \( s_U \) on a neighborhood of \( x \). Thus, \( s \) is continuous. Furthermore, for each \( U \in \mathscr{U} \), \( s|_U = s_U \). Hence, \( s \in \Gamma(O, \mathscr{T}) \), and the gluing property holds.

Therefore, \(\textsc{Ps}(\mathscr{T})\) is a sheaf of residuated lattices.
\end{proof}


The next result constructs a morphism of presheaves of residuated lattices arising from a morphism of \'{e}tal\'{e} spaces of residuated lattices.

\begin{proposition}\label{telpreshsecmo}
Let \( h \colon \mathscr{T} \overset{\pi}{\downarrow} \mathscr{B} \to \mathscr{S} \overset{\phi}{\downarrow} \mathscr{B} \) be a morphism of \'{e}tal\'{e} spaces of residuated lattices. For each open set \( U \subseteq \mathscr{B} \), define
\[
\textsc{Ps}(h)_U \colon \Gamma(U, \mathscr{T}) \to \Gamma(U, \mathscr{S}), \quad \sigma \mapsto h \circ \sigma.
\]
Then the family \( \textsc{Ps}(h) = \{ \textsc{Ps}(h)_U \}_{U \in \mathscr{O}(\mathscr{B})} \) defines a morphism of presheaves from \( \textsc{Ps}(\mathscr{T}) \) to \( \textsc{Ps}(\mathscr{S}) \).
\end{proposition}

\begin{proof}
For each open set \( U \subseteq \mathscr{B} \), the map \( \textsc{Ps}(h)_U \) is a morphism of residuated lattices because \( h \) is continuous and stalkwise a morphism of residuated lattices.

Now let \( U \subseteq V \subseteq \mathscr{B} \) and take \( \sigma \in \Gamma(V, \mathscr{T}) \). Then:
\[
\begin{aligned}
\textsc{Ps}(\mathscr{S})_{V,U} \big( \textsc{Ps}(h)_V(\sigma) \big) &= \textsc{Ps}(\mathscr{S})_{V,U}(h \circ \sigma) \\
&= (h \circ \sigma)|_U = h \circ (\sigma|_U) \\
&= \textsc{Ps}(h)_U(\sigma|_U) = \textsc{Ps}(h)_U \big( \textsc{Ps}(\mathscr{T})_{V,U}(\sigma) \big).
\end{aligned}
\]
Thus, the following diagram commutes:
\[
\xymatrix{
\Gamma(V,\mathscr{T}) \ar@{->}[dd]|-{{\textsc{Ps}(\mathscr{T})_{V,U}}} \ar@{->}[rr]|-{\textsc{Ps}(h)_{V}} &  & \Gamma(V,\mathscr{S}) \ar@{->}[dd]|-{{\textsc{Ps}(\mathscr{S})_{V,U}}} \\
 &  &  \\
\Gamma(U,\mathscr{T}) \ar@{->}[rr]|-{\textsc{Ps}(h)_{U}} &  & \Gamma(U,\mathscr{S})
}
\]
Hence, \( \textsc{Ps}(h) \) is a natural transformation of presheaves.
\end{proof}


\begin{theorem}\label{etalesheafunc}
For any topological space \( \mathscr{B} \), the assignment
\[
\textsc{Ps} \colon \mathbf{RL\text{-}Etale}(\mathscr{B}) \to \mathbf{RL\text{-}Sheaf}(\mathscr{B}),
\]
which sends each \'{e}tal\'{e} space of residuated lattices to its presheaf of sections, and each morphism \( h \) to the natural transformation \( \textsc{Ps}(h) \), is a functor.
\end{theorem}

\begin{proof}
This follows immediately from Theorem~\ref{etaleshea}, which ensures that \( \textsc{Ps}(\mathscr{T}) \) is indeed a sheaf, and from Proposition~\ref{telpreshsecmo}, which shows that \( \textsc{Ps}(h) \) defines a natural transformation between sheaves. Functoriality follows from the fact that identity maps and composition of morphisms are preserved under post-composition.
\end{proof}

\section{The Sheafification of Presheaves of Residuated Lattices}\label{sec5}

This section describes the sheafification process, which transforms a presheaf of residuated lattices into a sheaf of residuated lattices. Although this approach is often attributed to Michel Lazard \citep{cartan1950faisceaux}, we adopt the description given by Serre \citep{serre1955faisceaux}, which involves constructing a functor from \(\mathbf{RL\text{-}PreSheaf}(\mathscr{B})\) to \(\mathbf{RL\text{-}Etale}(\mathscr{B})\), for a fixed topological space \(\mathscr{B}\).

Recall from \cite[p.~11]{gratzer2008universal} that a poset \(\mathscr{A}\) is called \textit{directed} if every finite subset of \(\mathscr{A}\) has an upper bound in \(\mathscr{A}\).


\begin{examples}\label{direcposeexa}
Let \(\mathscr{B}\) be a topological space, and fix a subset \(X \subseteq \mathscr{B}\). The set of neighbourhoods of \(X\), denoted \(N(X)\), forms a directed poset under reverse inclusion.
\end{examples}


\begin{definition}\label{directsys}
Let \(I\) be a directed poset, \(\mathscr{A} = \{ \mathfrak{A}_i \}_{i \in I}\) a family of residuated lattices, and \(\mathscr{H} = \{ h_{ij} \colon \mathfrak{A}_i \to \mathfrak{A}_j \mid i \leq j \}\) a family of morphisms of residuated lattices which is called transition morphisms. The triple \((\mathscr{A}, \mathscr{H}, I)\) is called a \textit{direct system} of residuated lattices if satisfies the following assertions:
\begin{enumerate}
  \item [$(Id)$ \namedlabel{id}{$(Id)$}] The identity condition, i. e. \(h_{ii} = 1_{\mathfrak{A}_i}\) for every \(i \in I\);
  \item [$(Co)$ \namedlabel{co}{$(Co)$}] The compatibility condition, i. e. \(h_{ik} = h_{jk} \circ h_{ij}\) for all \(i \leq j \leq k\).
\end{enumerate}

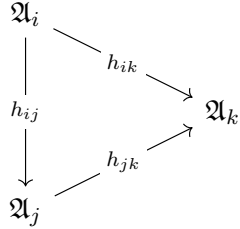
\begin{figure}[h!]
\centering
\begin{tikzcd}
\mathfrak{A}_i \arrow[rrd, "h_{ik}" description] \arrow[dd, "h_{ij}" description] &  &                  \\
                                                                                   &  & \mathfrak{A}_k \\
\mathfrak{A}_j \arrow[rru, "h_{jk}" description]                                  &  &
\end{tikzcd}
\end{figure}
\end{definition}


\begin{examples}\label{preshexa}
Let \(\mathscr{F}\) be a presheaf of residuated lattices over \(\mathscr{B}\), and fix a subset \(X \subseteq \mathscr{B}\). Define \(\mathscr{A} = \{ \mathscr{F}(U) \}_{U \in N(X)}\) and \(\mathscr{H} = \{ \mathscr{F}_{VU} \mid U \subseteq V \}\), where \(U, V \in N(X)\). Then the triple
\[
\mathscr{D}_X^{\mathscr{F}} := (\mathscr{A}, \mathscr{H}, N(X))
\]
is a direct system of residuated lattices. In the special case where \(X = \{x\}\) for some \(x \in \mathscr{B}\), we denote the system by \(\mathscr{D}_x^{\mathscr{F}}\).
\end{examples}


\begin{definition}\label{directlimitdef}
Let \(\mathscr{D} = (\mathscr{A}, \mathscr{H}, I)\) be a direct system of residuated lattices, where \(\mathscr{A} = \{ \mathfrak{A}_i \}_{i \in I}\), and \(\mathscr{H} = \{ h_{ij} \colon \mathfrak{A}_i \to \mathfrak{A}_j \mid i \leq j \}\). Let \(\mathfrak{A}\) be a residuated lattice and let \(\psi = \{ \psi_i \colon \mathfrak{A}_i \to \mathfrak{A} \}_{i \in I}\) be a family of morphisms of residuated lattices. Then the pair \((\mathfrak{A}, \psi)\) is called a \textit{direct limit} of the direct system \(\mathscr{D}\) if the following conditions hold:
\begin{enumerate}
  \item [$(Co)$ \namedlabel{cod}{$(Co)$}] The compatibility condition, i. e. $\psi_{i}=\psi_{j}h_{ij}$ whenever $i\leq j$;
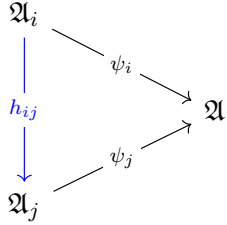
\begin{figure}[h!]
\centering
\begin{tikzcd}
\mathfrak{A}_{i} \arrow[rrd, "\psi_{i}" description] \arrow[dd, "h_{ij}" description,blue] &  &              \\
                                                                                      &  & \mathfrak{A} \\
\mathfrak{A}_{j} \arrow[rru, "\psi_{j}" description]                                  &  &
\end{tikzcd}
\end{figure}
\item [$(Un)$ \namedlabel{un}{$(Un)$}] The universal property, i. e.  for any residuated lattice $\mathfrak{B}$ and any family of morphisms of residuated lattices $\{\varphi_{i}:\mathfrak{A}_{i}\longrightarrow \mathfrak{B}\}_{i\in I}$ such that $\varphi_{i}=\varphi_{j}h_{ij}$ whenever $i\leq j$, there exists a unique morphisms of residuated lattices $h:\mathfrak{A}\longrightarrow \mathfrak{B}$ such that $\varphi_{i}=h\psi_{i}$, for any $i\in I$.
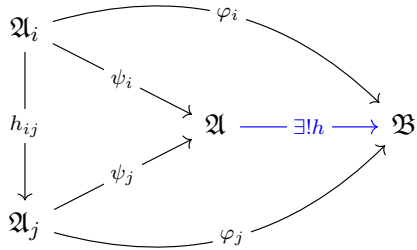
\begin{figure}[h!]
\centering
\begin{tikzcd}
\mathfrak{A}_{i} \arrow[rrd, "\psi_{i}" description] \arrow[dd, "h_{ij}" description] \arrow[rrrrd, "\varphi_{i}" description, bend left] &  &                                             &  &              \\
&  & \mathfrak{A}_{} \arrow[rr, "\exists!h" description,blue] &  & \mathfrak{B} \\
\mathfrak{A}_{j} \arrow[rru, "\psi_{j}" description] \arrow[rrrru, "\varphi_{j}" description, bend right]                                 &  &                                             &  &
\end{tikzcd}
\end{figure}
\end{enumerate}
\end{definition}


The following result establishes the existence and uniqueness (up to isomorphism) of the direct limit for a direct system of residuated lattices.

\begin{proposition}\label{uniqdirlim}
Let $\mathscr{D} = (\mathscr{A}, \mathscr{H}, I)$ be a direct system of residuated lattices. Then a direct limit of $\mathscr{D}$ exists and is unique up to isomorphism.
\end{proposition}

\begin{proof}
Define a relation $\equiv$ on the disjoint union $\coprod \mathscr{A}$ by:
\[
x \equiv y \quad \text{if and only if} \quad h_{ik}(x) = h_{jk}(y)
\]
for some indices $i, j \leq k$, where $x \in A_i$ and $y \in A_j$.

It is straightforward to verify that $\equiv$ is an equivalence relation. Let $\mathscr{D}[x]$ denote the equivalence class of $x$, and define the set
\[
\varinjlim \mathscr{D} := \left\{ \mathscr{D}[x] \,\middle|\, x \in \coprod \mathscr{A} \right\}.
\]

For any $\mathscr{D}[x], \mathscr{D}[y] \in \varinjlim \mathscr{D}$ and any binary fundamental operation $\diamond$ of residuated lattices, define:
\[
\mathscr{D}[x] \diamond \mathscr{D}[y] := \mathscr{D}[h_{ik}(x) \diamond_k h_{jk}(y)],
\]
where $x \in A_i$, $y \in A_j$, and $i, j \leq k$. It is obvious that the choice of $k$ does not matter due to the compatibility of the system.

Define constants $\tilde{0}, \tilde{1} \in \varinjlim \mathscr{D}$ by
\[
\tilde{0} := \mathscr{D}[0_i], \quad \tilde{1} := \mathscr{D}[1_i]
\]
for any $i \in I$.

One verifies that $\left( \varinjlim \mathscr{D}; \vee, \wedge, \odot, \to, \tilde{0}, \tilde{1} \right)$ forms a residuated lattice. Define canonical morphisms $\psi_i: A_i \rightarrow \varinjlim \mathscr{D}$ by $\psi_i(x) := \mathscr{D}[x]$ for all $x \in A_i$. Then $(\varinjlim \mathscr{D}, \{ \psi_i \}_{i \in I})$ is a direct limit of $\mathscr{D}$.

The uniqueness of the direct limit, up to isomorphism, follows from the universal property and is omitted.
\end{proof}


\begin{lemma}\label{lemmohem}
Let $(\mathscr{A}, \mathscr{H}, I)$ be a direct system of residuated lattices. If $i, j \in I$ with $i \leq j$ and $x \in \mathfrak{A}_i$, then $x \equiv h_{ij}(x)$.
\end{lemma}

\begin{proof}
This follows immediately by taking $k = j$ in the definition of the relation $\equiv$.
\end{proof}


\begin{definition}
Let $\mathscr{F}$ be a presheaf of residuated lattices on a topological space $\mathscr{B}$, and let $X \subseteq \mathscr{B}$. The direct limit $\varinjlim \mathscr{D}_X^{\mathscr{F}}$ is denoted by $\mathscr{F}_X$.

For any neighborhood $U$ of $X$ in $\mathscr{B}$ and any section $s \in \mathscr{F}(U)$, the equivalence class $\mathscr{D}_X^{\mathscr{F}}[s]$ is denoted by $\mathscr{F}_{U,X}(s)$, or simply $[s]_X$ when the domain is clear.

For $x \in \mathscr{B}$, the residuated lattice $\mathscr{F}_{\{x\}}$ is denoted by $\mathscr{F}_x$ and is called the \emph{stalk} of $\mathscr{F}$ at $x$. The class $\mathscr{F}_{U,\{x\}}(s)$ is denoted by $\mathscr{F}_{U,x}(s)$, or simply $[s]_x$, and is referred to as the \emph{germ} of $s$ at $x$.
\end{definition}

Intuitively, the germ of a section at a point in the base space captures the local behavior of that section within an arbitrarily small neighborhood around the point, while disregarding how the section behaves elsewhere. More precisely, the germ encodes the infinitesimal information of a section at a given point, abstracting away any global variation.

Two sections are considered equivalent at a point if they agree on some open neighborhood containing that point. The germ of a section is thus the equivalence class of all such sections under this relation. In this way, the germ represents the local identity of a section at a point, independent of its behavior outside any neighborhood of that point.

In contrast, the stalk of a presheaf of residuated lattices at a point in the base space is the collection of all germs at that point. Intuitively, the stalk encapsulates all possible local behaviors of sections near the point. One can think of the stalk as the result of ``zooming in'' infinitely close to a point: it records the local values that sections can attain near the point, while ignoring distinctions that arise farther away.


\begin{example}\label{stalkgermex}
Consider the presheaf \( \mathscr{F} \) described in Example~\ref{exsirsheanot}. The stalk \( \mathscr{F}_x \) at the point \( x \) is formed by taking equivalence classes of sections defined on neighborhoods of \( x \), where two sections are considered equivalent if they agree on the value of \( \alpha \) when restricted to \( \{x\} \). Since every section over \( \{x\} \) is required to satisfy \( \nu(\beta) = 0 \), the germ at \( x \) is completely determined by the value of \( \alpha \). 

In particular, for any section \( s \in \mathscr{F}(\{x\}) \), the germ of \( s \) at \( x \) is given by
\[
[s]_x = \{ t \in \mathscr{F}(\{x,y\}) \mid t(\alpha) = s(\alpha) \}.
\]
Therefore, the stalk at \( x \) is
\[
\mathscr{F}_x = \{ [s]_x \mid s \in \mathscr{F}(\{x\}) \} \cong [0,1],
\]
with the structure of a residuated lattice given by the \L{}ukasiewicz operations. Each germ at \( x \) corresponds uniquely to a value \( \nu(\alpha) \in [0,1] \), and the value of \( \beta \) is always zero by construction.

On the other hand, the stalk \( \mathscr{F}_y \) at the point \( y \) is determined entirely by the values of sections in \( \mathscr{F}(\{x, y\}) \), since \( \{x, y\} \) is the only open set that contains \( y \). Each such section assigns values to both \( \alpha \) and \( \beta \), and since no additional restriction applies at \( y \), every section forms a singleton germ. That is,
\[
[s]_y = \{ s \},
\]
and the stalk is given by
\[
\mathscr{F}_y \cong [0,1] \times [0,1],
\]
with pointwise (componentwise) \L{}ukasiewicz operations. Each germ at \( y \) thus corresponds to a pair \( (\nu(\alpha), \nu(\beta)) \in [0,1]^2 \).

\end{example}


\begin{lemma}\label{kheilimohem}
Let $\mathscr{F}$ be a presheaf of residuated lattices on $\mathscr{B}$, and let $U \subseteq V$ be neighborhoods of $X$ in $\mathscr{B}$. Then, for any $s \in \mathscr{F}(V)$,
\[
\mathscr{F}_{V,X}(s) = \mathscr{F}_{U,X}(s|_U).
\]
\end{lemma}

\begin{proof}
This is an immediate consequence of Lemma~\ref{lemmohem}.
\end{proof}


\begin{proposition}\label{sectgermequ}
  Let $\mathscr{F}$ be a presheaf of residuated lattices over $\mathscr{B}$ satisfying condition \ref{s} , $U$ an open subset of $\mathscr{B}$ and $s,t\in \mathscr{F}(U)$. The following assertions are equivalent:
  \begin{enumerate}
  \item  [$(1)$ \namedlabel{sectgermequ1}{$(1)$}] $s=t$;
  \item  [$(2)$ \namedlabel{sectgermequ2}{$(2)$}] $[s]_{x}=[t]_{x}$, for all $x\in U$.
\end{enumerate}
\end{proposition}
\begin{proof}
  \item [\ref{sectgermequ1}$\Rightarrow$\ref{sectgermequ2}:] It is evident.
  \item [\ref{sectgermequ2}$\Rightarrow$\ref{sectgermequ1}:] For any $x\in U$ there exists a neighborhood $U_{x}$ of $x$ contained in $U$ such that $s|_{U_{x}}=t|_{U_{x}}$. Since $\{U_{x}\}_{x\in U}$ is an open covering of $U$, the result is obtained.
\end{proof}


The rest of this section is devoted to constructing a functor from the category $\textbf{RL-PreSheaf}(\mathscr{B})$ to the category $\textbf{RL-Etale}(\mathscr{B})$.

\begin{proposition}\label{shefifipro}
Let \( \mathscr{F} \) be a presheaf of residuated lattices on a topological space \( \mathscr{B} \). For any open set \( U \subseteq \mathscr{B} \) and any section \( s \in \mathscr{F}(U) \), define the \emph{germ map} associated to \( s \) as the function
\[
s_U : U \to \mathrm{Et}(\mathscr{F}), \quad s_U(x) := [s]_x,
\]
where \( \mathrm{Et}(\mathscr{F}) := \coprod_{x \in \mathscr{B}} \mathscr{F}_x \) denotes the \'{e}tal\'{e} space of \( \mathscr{F} \), constructed as the disjoint union of all stalks.

The collection of subsets
\[
\left\{ \operatorname{Im}(s_U) \;\middle|\; U \in \mathscr{O}(\mathscr{B}),\ s \in \mathscr{F}(U) \right\}
\]
forms a basis for the \emph{final topology} on \( \mathrm{Et}(\mathscr{F}) \) induced by the family of maps \( \{ s_U \}_{U \in \mathscr{O}(\mathscr{B})} \).
\end{proposition}

\begin{proof}
Let $U,V$ be open sets in $\mathscr{B}$, $s\in \mathscr{F}(U)$ and $t\in \mathscr{F}(V)$. With a little bit of effort is inferred that $b\in t_{V}^{\leftarrow}(Im(s_{U}))$ if and only if there exists a neighborhood $W_{b}$ of $b$ such that it is contained in $U\cap V$ and $s|_{W_{b}}=t|_{W_{b}}$. One can see that $t_{V}^{\leftarrow}(Im(s_{U}))=\bigcup_{b\in t_{V}^{\leftarrow}(Im(s_{U}))}W_{b}$. Thus, $\operatorname{Im}(s_U)$ is open in the final topology.

Now, let $O \subseteq Et(\mathscr{F})$ be open in the final topology, and suppose $\mathscr{F}_{U,x}(s) \in O$. Then $x \in s_U^{-1}(O)$, and since this set is open in $U$, there exists an open neighborhood $W \subseteq U$ of $x$ such that
\[
\mathscr{F}_{U,x}(s) = \mathscr{F}_{W,x}(s|_W) \in \operatorname{Im}(s|_W)_W \subseteq O.
\]
Therefore, the sets $\operatorname{Im}(s_U)$ form a basis for the final topology.
\end{proof}


Let $\mathscr{F}$ be a presheaf of residuated lattices over the topological space $\mathscr{B}$. We equip the set $Et(\mathscr{F})$ with the topology generated by the basis 
\[
\{ \operatorname{Im}(s_U) \mid U \in \mathscr{O}(\mathscr{B}),\ s \in \mathscr{F}(U) \},
\]
referred to as the \emph{\'{e}tal\'{e} topology}. We denote by 
\[
\pi_{\mathscr{F}}: Et(\mathscr{F}) \longrightarrow \mathscr{B}
\]
the canonical projection given by $\pi_{\mathscr{F}}([s]_x) = x$.

\begin{lemma}\label{preshetaov}
Let $\mathscr{F}$ be a presheaf of residuated lattices over $\mathscr{B}$. Then the pair $(Et(\mathscr{F}), \pi_{\mathscr{F}})$ is an \'{e}tal\'{e} space over $\mathscr{B}$.
\end{lemma}

\begin{proof}
To prove that $\pi_{\mathscr{F}}$ is continuous, observe that for any open set $U \in \mathscr{O}(\mathscr{B})$, the preimage under $\pi_{\mathscr{F}}$ is given by
\[
\pi_{\mathscr{F}}^{-1}(U) = \bigcup \left\{ \operatorname{Im}(s_V) \mid V \in \mathscr{O}(\mathscr{B}),\ V \subseteq U,\ s \in \mathscr{F}(V) \right\}.
\]
Hence, $\pi_{\mathscr{F}}$ is continuous. 

Next, let $U$ be an open subset of $\mathscr{B}$ and let $s \in \mathscr{F}(U)$. Then 
\[
\pi_{\mathscr{F}}(\operatorname{Im}(s_U)) = U,
\]
which shows that $\pi_{\mathscr{F}}$ is an open map. 

Now consider an element $[s]_{x} \in Et(\mathscr{F})$. Since $\pi_{\mathscr{F}}|_{\operatorname{Im}(s_U)}$ is injective, it follows that $\pi_{\mathscr{F}}$ is locally injective. By Theorem~\ref{lhcol}, this implies that $\pi_{\mathscr{F}}$ is a local homeomorphism.
\end{proof}


\begin{theorem}\label{shefifipro1}
Let $\mathscr{F}$ be a presheaf of residuated lattices over $\mathscr{B}$. The \'{e}tal\'{e} space $Et(\mathscr{F})\overset{\pi_{\mathscr{F}}}{\downarrow}\mathscr{B}$ is an \'{e}tal\'{e} space of residuated lattices.
\end{theorem}
\begin{proof}
  Applying Lemma \ref{preshetaov}, the pair \( (Et(\mathscr{F}), \pi_{\mathscr{F}}) \) is an \'{e}tal\'{e} space over \( \mathscr{B} \). For each \( x \in \mathscr{B} \), we have \( \textsc{Et}(\mathscr{F})_x = \mathscr{F}_x \), and thus the stalks of \( Et(\mathscr{F}) \) are residuated lattices. Assume that \( \diamond \) is a binary fundamental operation. Considering
\[
Et(\mathscr{F})^{(2)} = \{ (\mathscr{F}_{U,x}(s), \mathscr{F}_{V,x}(t)) \mid x \in \mathscr{B} \},
\]
we have:
\[
\overset{Et(\mathscr{F})}{\diamond}: (\mathscr{F}_{U,x}(s), \mathscr{F}_{V,x}(t)) \rightsquigarrow \mathscr{F}_{W,x}\left( \overset{\mathscr{F}(W)}{\diamond}(s|_W, t|_W) \right),
\]
where \( W \) is an appropriate neighborhood of \( x \). Let \( \mathscr{F}_{U,x}(s), \mathscr{F}_{V,x}(t) \in Et(\mathscr{F}) \). Thus, there exists a neighborhood \( W \) of \( x \), contained in \( U \) and \( V \), such that
\[
\overset{Et(\mathscr{F})}{\diamond}(\mathscr{F}_{U,x}(s), \mathscr{F}_{V,x}(t)) = \overset{\mathscr{F}(W)}{\diamond}(s|_W, t|_W).
\]
Assume \( \mathscr{O} \) is a neighborhood of \( \mathscr{F}_{W,x}\left( \overset{\mathscr{F}(W)}{\diamond}(s|_W, t|_W) \right) \). Then, there exists an open subset \( H \subseteq \mathscr{B} \) and a section \( r \in \mathscr{F}(H) \) such that \( \text{Im}(r_H) \) is a basic neighborhood of \( \mathscr{F}_{W,x}\left( \overset{\mathscr{F}(W)}{\diamond}(s|_W, t|_W) \right) \), contained in \( \mathscr{O} \). This implies
\[
\mathscr{F}_{W,x}\left( \overset{\mathscr{F}(W)}{\diamond}(s|_W, t|_W) \right) = \mathscr{F}_{H,x}(r).
\]
Thus, there exists a neighborhood \( L \) of \( x \), contained in both \( W \) and \( H \), such that
\[
\overset{\mathscr{F}(W)}{\diamond}(s|_W, t|_W) |_L = r|_L.
\]
By Lemma \ref{kheilimohem}, the pullback \( \text{Im}((s|_L)_L) \times_{\mathscr{B}} \text{Im}((t|_L)_L) \) forms an open neighborhood of \( (\mathscr{F}_{U,x}(s), \mathscr{F}_{V,x}(t)) \). Let \( (\mathscr{F}_{L,y}(s|_L), \mathscr{F}_{L,y}(t|_L)) \) be an element of
\[
\text{Im}((s|_L)_L) \times_{\mathscr{B}} \text{Im}((t|_L)_L).
\]
We then have the following formulas:
\[
\begin{array}{ll}
\overset{Et(\mathscr{F})}{\diamond}(\mathscr{F}_{L,y}(s|_L), \mathscr{F}_{L,y}(t|_L)) & = \mathscr{F}_{L,y}\left( \overset{\mathscr{F}(L)}{\diamond}(s|_L, t|_L) \right) \\
& = \mathscr{F}_{L,y}\left( (\overset{\mathscr{F}(W)}{\diamond}(s|_W, t|_W)) |_L \right) \\
& = \mathscr{F}_{L,y}(r|_L) \\
& = \mathscr{F}_{H,y}(r).
\end{array}
\]
This shows that
\[
\overset{Et(\mathscr{F})}{\diamond}\left( \text{Im}((s|_L)_L) \times_{\mathscr{B}} \text{Im}((t|_L)_L) \right) \subseteq \text{Im}(r_H).
\]
Thus, the operation \( \overset{Et(\mathscr{F})}{\diamond} \) is continuous.

For the case where \( \diamond \) is a nullary fundamental operation, we have \( \overset{f^*\mathscr{S}}{\diamond} \) as the section \( \diamond_B \), and it follows that it is continuous.

\end{proof}


\begin{definition}
  Let \( \mathscr{F} \) be a presheaf of residuated lattices over \( \mathscr{B} \). In the sequel, the \'{e}tal\'{e} space of residuated lattices
\[
Et(\mathscr{F}) \overset{\pi_{\mathscr{F}}}{\downarrow} \mathscr{B}
\]
is called the \textit{stalk space of \( \mathscr{F} \)}.
\end{definition}


The following example clarifies the above discussion by explicitly constructing the stalk space \( \mathrm{Et}(\mathscr{F}) \) of the presheaf of residuated lattices \( \mathscr{F} \) over the Sierpi\'{n}ski space, as described in Example~\ref{exsirsheanot}.

\begin{example}
Consider the presheaf \( \mathscr{F} \) defined in Example~\ref{exsirsheanot}. According to Example~\ref{stalkgermex}, the stalk at the point \( x \), denoted \( \mathscr{F}_x \), is isomorphic to the unit interval \( [0,1] \), consisting of germs determined solely by the value of \( \alpha \), since every section over an open set containing \( x \) satisfies \( \nu(\beta) = 0 \). In contrast, the stalk at the point \( y \), denoted \( \mathscr{F}_y \), is isomorphic to the Cartesian product \( [0,1] \times [0,1] \), as both \( \alpha \) and \( \beta \) are observable and unconstrained at \( y \), through sections defined over \( \{x, y\} \).

The \'{e}tal\'{e} space of \( \mathscr{F} \) is the disjoint union
\[
\mathrm{Et}(\mathscr{F}) = \mathscr{F}_x \coprod \mathscr{F}_y \cong [0,1] \coprod [0,1]^2,
\]
equipped with the \emph{\'{e}tal\'{e} topology}, as described in Proposition~\ref{shefifipro}. This topology is generated by the images of germ maps associated with local sections of the presheaf.

Specifically, for a section \( s \in \mathscr{F}(\{x\}) \), the germ map \( s_{\{x\}} \colon \{x\} \to \mathrm{Et}(\mathscr{F}) \) has image
\[
\operatorname{Im}(s_{\{x\}}) = \{ [s]_x \} \subset \mathscr{F}_x,
\]
corresponding to the singleton point \( s(\alpha) \in [0,1] \), since \( \beta \) is forced to be zero and plays no role in the germ at \( x \).

For a section \( s \in \mathscr{F}(\{x, y\}) \), the germ map \( s_{\{x,y\}} \colon \{x, y\} \to \mathrm{Et}(\mathscr{F}) \) has image
\[
\operatorname{Im}(s_{\{x,y\}}) = \{ [s]_x,\ [s]_y \},
\]
where \( [s]_x \in \mathscr{F}_x \) corresponds to the germ determined by \( \nu(\alpha) \), with \( \nu(\beta) = 0 \) implicitly enforced, and \( [s]_y \in \mathscr{F}_y \) corresponds to the full pair \( (\nu(\alpha), \nu(\beta)) \in [0,1]^2 \). Thus, the image of this germ map links a point in \( \mathscr{F}_x \) with a corresponding point in \( \mathscr{F}_y \), forming a basic open set in the \'{e}tal\'{e} topology.

The \'{e}tal\'{e} projection
\[
\pi: \mathrm{Et}(\mathscr{F}) \longrightarrow \mathscr{B}
\]
maps each germ \( [s]_z \) to the point \( z \in \mathscr{B} \) at which it is taken; that is, \( \pi([s]_x) = x \) and \( \pi([s]_y) = y \). By Theorem~\ref{shefifipro1}, this projection equips \( \mathrm{Et}(\mathscr{F}) \) with the structure of an \'{e}tal\'{e} space of residuated lattices.

This construction illustrates how the local algebraic data encoded in the presheaf \( \mathscr{F} \) is spatially organized: each point in the base space \( \mathscr{B} \) is associated with a stalk (a local algebra), and the \'{e}tal\'{e} topology coherently glues these stalks via shared local sections. As a result, the operations of the residuated lattice structure remain continuous when viewed through the \'{e}tal\'{e} space, and \( \pi \) becomes a local homeomorphism.

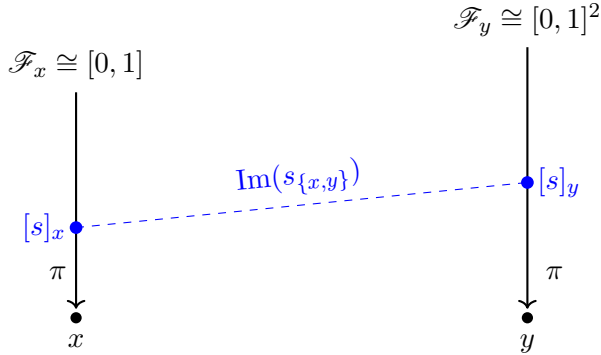
\begin{center}
\begin{tikzpicture}[scale=1.2]
  \node[fill=black, circle, inner sep=1.5pt, label=below:$x$] (x) at (0,0) {};
  \node[fill=black, circle, inner sep=1.5pt, label=below:$y$] (y) at (5,0) {};
  
  \draw[->, thick] (0,2.5) -- (0,0.1);
  \draw[->, thick] (5,3) -- (5,0.1);
  \node at (-0.2,0.5) {$\pi$};
  \node at (5.3,0.5) {$\pi$};

  \node at (0,2.8) {$\mathscr{F}_x \cong [0,1]$};
  \node at (5,3.3) {$\mathscr{F}_y \cong [0,1]^2$};

  \fill[blue] (0,1) circle (2pt) node[left] {$[s]_x$};
  \fill[blue] (5,1.5) circle (2pt) node[right] {$[s]_y$};

  \draw[dashed, blue] (0,1) -- (5,1.5) node[midway, above, sloped] {$\operatorname{Im}(s_{\{x,y\}})$};
\end{tikzpicture}
\end{center}
\end{example}


The following result constructs a morphism between \'{e}tal\'{e} spaces of residuated lattices, induced by a morphism of presheaves of residuated lattices.

\begin{theorem}\label{etpresheafmor}
Let \( \mathscr{F} \) and \( \mathscr{G} \) be presheaves of residuated lattices over the topological space \( \mathscr{B} \), and let \( \phi: \mathscr{F} \longrightarrow \mathscr{G} \) be a morphism of presheaves of residuated lattices over \( \mathscr{B} \). Then the map
\[
\textsc{Et}(\phi): \textsc{Et}(\mathscr{F}) \longrightarrow \textsc{Et}(\mathscr{G}),
\]
defined by
\[
\mathscr{F}_{U,x}(s) \mapsto \mathscr{G}_{U,x}(\phi(U)(s)),
\]
is a morphism of \'{e}tal\'{e} spaces of residuated lattices.
\end{theorem}

\begin{proof}
Let \( \mathscr{F}_{U,x}(s) = \mathscr{F}_{V,y}(t) \). Then clearly \( x = y \). Hence, there exists an open neighborhood \( W \subseteq U \cap V \) of \( x \) such that 
\[
\mathscr{F}_{UW}(s) = \mathscr{F}_{VW}(t).
\]
Applying the functoriality of \( \phi \), we obtain:
\[
\mathscr{G}_{UW}(\phi(U)(s)) = \phi(W)(\mathscr{F}_{UW}(s)) = \phi(W)(\mathscr{F}_{VW}(t)) = \mathscr{G}_{VW}(\phi(V)(t)).
\]
Therefore,
\[
\textsc{Et}(\phi)(\mathscr{F}_{U,x}(s)) = \mathscr{G}_{U,x}(\phi(U)(s)) = \mathscr{G}_{V,x}(\phi(V)(t)) = \textsc{Et}(\phi)(\mathscr{F}_{V,x}(t)),
\]
so \( \textsc{Et}(\phi) \) is well-defined on germs.

Let \( U \subseteq \mathscr{B} \) be open and let \( s \in \mathscr{F}(U) \). Define \( t = \phi(U)(s) \in \mathscr{G}(U) \). Then
\[
\textsc{Et}(\phi)(\operatorname{Im} s_U) = \operatorname{Im}(t_U),
\]
so \( \textsc{Et}(\phi) \) maps basic open sets to basic open sets, and is therefore an open map. It is also straightforward to verify that
\[
\pi_{\mathscr{G}} \circ \textsc{Et}(\phi) = \pi_{\mathscr{F}},
\]
and thus \( \textsc{Et}(\phi) \) is continuous by Proposition~\ref{conopnloc}.

Next, we show that for each \( x \in \mathscr{B} \), the map
\[
\textsc{Et}(\phi)_x: \textsc{Et}(\mathscr{F})_x \longrightarrow \textsc{Et}(\mathscr{G})_x
\]
is a morphism of residuated lattices. Let \( \mathscr{F}_{U,x}(s), \mathscr{F}_{V,x}(t) \in \mathscr{F}_x \). Then there exists a neighborhood \( W \subseteq U \cap V \) of \( x \) such that
\[
\mathscr{F}_{U,x}(s) = \mathscr{F}_{W,x}(s|_W), \quad \mathscr{F}_{V,x}(t) = \mathscr{F}_{W,x}(t|_W).
\]
Let \( \diamond \) be any binary operation in the signature of residuated lattices (e.g., \( \wedge, \vee, \cdot, \rightarrow \)). Then:
\[
\begin{aligned}
\textsc{Et}(\phi)_x(\mathscr{F}_{U,x}(s) \diamond \mathscr{F}_{V,x}(t)) 
&= \textsc{Et}(\phi)_x(\mathscr{F}_{W,x}(s|_W \diamond t|_W)) \\
&= \mathscr{G}_{W,x}(\phi(W)(s|_W \diamond t|_W)) \\
&= \mathscr{G}_{W,x}(\phi(W)(s|_W) \diamond \phi(W)(t|_W)) \\
&= \mathscr{G}_{W,x}(\phi(W)(s|_W)) \diamond \mathscr{G}_{W,x}(\phi(W)(t|_W)) \\
&= \textsc{Et}(\phi)_x(\mathscr{F}_{W,x}(s|_W)) \diamond \textsc{Et}(\phi)_x(\mathscr{F}_{W,x}(t|_W)) \\
&= \textsc{Et}(\phi)_x(\mathscr{F}_{U,x}(s)) \diamond \textsc{Et}(\phi)_x(\mathscr{F}_{V,x}(t)).
\end{aligned}
\]

Finally, to show that \( \textsc{Et}(\phi)_x \) preserves the nullary operations, let \( 0, 1 \in \mathscr{F}(U) \). Then by the property of morphisms of residuated lattices,
\[
\phi(U)(0) = 0_{\mathscr{G}(U)}, \quad \phi(U)(1) = 1_{\mathscr{G}(U)},
\]
and so
\[
\textsc{Et}(\phi)_x([0]_x) = [0]_x, \quad \textsc{Et}(\phi)_x([1]_x) = [1]_x.
\]
Therefore, \( \textsc{Et}(\phi)_x \) is a morphism of residuated lattices for each \( x \in \mathscr{B} \).
\end{proof}


\begin{theorem}\label{preshetalfunc}
For a topological space $\mathscr{B}$, the assignment
\[
\textsc{Et}: \mathbf{RL\text{-}PreSheaf}(\mathscr{B}) \longrightarrow \mathbf{RL\text{-}Etale}(\mathscr{B})
\]
defines a functor. Furthermore, if $\mathscr{G}$ satisfies condition~\ref{s}, then this functor is faithful on $\operatorname{Hom}(\mathscr{F}, \mathscr{G})$.
\end{theorem}

\begin{proof}
Let $\phi: \mathscr{F} \longrightarrow \mathscr{G}$ and $\psi: \mathscr{G} \longrightarrow \mathscr{H}$ be morphisms of presheaves of residuated lattices. The functoriality condition
\[
\textsc{Et}(\psi \circ \phi) = \textsc{Et}(\psi) \circ \textsc{Et}(\phi)
\]
is straightforward to verify. Moreover, the identity morphism is preserved: $\textsc{Et}(\mathrm{id}_{\mathscr{F}}) = \mathrm{id}_{\textsc{Et}(\mathscr{F})}$. Thus, by Theorems~\ref{shefifipro1} and~\ref{etpresheafmor}, it follows that
\[
\textsc{Et}: \mathbf{RL\text{-}PreSheaf}(\mathscr{B}) \longrightarrow \mathbf{RL\text{-}Etale}(\mathscr{B})
\]
is indeed a functor.

Now assume that $\mathscr{G}$ satisfies condition~\ref{s}, and let $\phi, \psi \in \operatorname{Hom}(\mathscr{F}, \mathscr{G})$ be such that $\textsc{Et}(\phi) = \textsc{Et}(\psi)$. We must show that $\phi = \psi$. 

Let $U$ be an open subset of $\mathscr{B}$ and let $s \in \mathscr{F}(U)$. For any $x \in U$, we have
\[
\textsc{Et}(\phi)([s]_x) = \textsc{Et}(\psi)([s]_x),
\]
which implies that
\[
[\phi(U)(s)]_x = [\psi(U)(s)]_x.
\]
By the definition of germ equivalence, there exists a neighborhood $U_x \subseteq U$ of $x$ such that
\[
\phi(U)(s)\big|_{U_x} = \psi(U)(s)\big|_{U_x}.
\]
The collection $\{U_x\}_{x \in U}$ forms an open cover of $U$, and by the uniqueness of sections over open covers (condition~\ref{s}), it follows that $\phi(U)(s) = \psi(U)(s)$. Hence, $\phi = \psi$, and the functor is faithful on $\operatorname{Hom}(\mathscr{F}, \mathscr{G})$.
\end{proof}


 \begin{definition}
Let $\mathscr{B}$ be a topological space. The composite functor
\[
\textsc{Ps} \circ \textsc{Et} : \mathbf{RL\text{-}PreSheaf}(\mathscr{B}) \longrightarrow \mathbf{RL\text{-}Sheaf}(\mathscr{B})
\]
is called the \emph{sheafification functor}. For any presheaf $\mathscr{F}$, the sheaf $\textsc{Ps}(\textsc{Et}(\mathscr{F}))$ is referred to as the \emph{sheafification of $\mathscr{F}$} and is denoted by $\mathscr{F}^+$.
\end{definition}


\begin{remark}\label{sheafifirema}
The sheafification functor \( \mathscr{F} \mapsto \mathscr{F}^+ \) can be intuitively understood as a \emph{completion process} that systematically repairs the gluing failures of a presheaf while preserving its underlying algebraic structure. Geometrically, this transformation converts local data—potentially incompatible across overlaps—into a globally coherent object through the following conceptual stages:

The presheaf \( \mathscr{F} \) is first decomposed into its germs (Theorem \ref{shefifipro1}), which represent the most refined local behaviors of sections at each point. This step highlights where the presheaf fails to satisfy the sheaf condition, i.e., where local compatibility breaks down.

A sheaf is then reconstructed by assembling these germs into a coherent structure, subject to the requirement that sections must agree on overlaps. The resulting sheaf \( \mathscr{F}^+ \) contains exactly those sections that can be patched together consistently from local data.

Algebraically, this process preserves the residuated lattice structure at each stalk while eliminating contradictions between incompatible local descriptions. As shown in Theorem~\ref{reflshbpshb}, the sheafification functor satisfies a universal property: it modifies a presheaf only to the extent necessary to satisfy the sheaf condition—no more, no less. Hence, \( \mathscr{F}^+ \) serves as the most efficient and canonical extension of \( \mathscr{F} \) into a sheaf, retaining all compatible local information and supplementing only what is required for global consistency.

In the context of fuzzy logic and substructural logics modeled by residuated lattices, sheafification provides a principled method to:

\begin{itemize}
\item Repair local truth assignments that fail to respect logical consequence on overlapping domains,
\item Extend partial valuations into coherent global truth-value assignments, preserving connectives such as \( \wedge, \vee, \odot, \to \),
\item Resolve semantic inconsistencies by enforcing coherence conditions across contexts.
\end{itemize}

For example, in graded modal logic, sheafification ensures that presheaves assigning inconsistent degrees of truth to related propositions (e.g., across different possible worlds) are corrected to satisfy the \emph{heredity condition} required by \( K \)-frames. This process is foundational in the development of sheaf semantics for non-classical logics.
\end{remark}


\begin{example}
Let \( \mathscr{F} \) be the presheaf defined in Example~\ref{exsirsheanot}, and let \( \mathrm{Et}(\mathscr{F}) \to \mathscr{B} \) denote its \'{e}tal\'{e} space, where \( \mathscr{B} = \{x, y\} \) with topology \( \mathscr{O}(\mathscr{B}) = \{\varnothing, \{x\}, \{x,y\}\} \). We now construct the presheaf of sections \( \textsc{Ps}(\mathrm{Et}(\mathscr{F})) \), as in Proposition~\ref{contavafunctshe}.

To each open set \( U \subseteq \mathscr{B} \), we assign the set \( \Gamma(U, \mathrm{Et}(\mathscr{F})) \) of continuous sections \( \sigma: U \to \mathrm{Et}(\mathscr{F}) \) satisfying \( \pi \circ \sigma = \mathrm{id}_U \).

We compute \( \textsc{Ps}(\mathrm{Et}(\mathscr{F})) \) explicitly:

\smallskip
\noindent
For \( U = \varnothing \) there are no points to map, so the only section is the empty function:
\[
\textsc{Ps}(\mathrm{Et}(\mathscr{F}))(\varnothing) = \{\varnothing\},
\]
which corresponds to the one-element (trivial) residuated lattice.

\smallskip
\noindent
For \( U = \{x\} \) a section selects a germ \( [s]_x \in \mathscr{F}_x \cong [0,1] \), determined solely by the value \( \nu(\alpha) \), since \( \nu(\beta) = 0 \) is fixed by the presheaf construction. Thus,
\[
\textsc{Ps}(\mathrm{Et}(\mathscr{F}))(\{x\}) \cong [0,1],
\]
with the $\L$ukasiewicz residuated lattice structure.

\smallskip
\noindent
For \( U = \{x,y\} \) a section corresponds to a single section \( s \in \mathscr{F}(\{x,y\}) \), and its germ map yields
\[
\sigma(x) = [s]_x \in \mathscr{F}_x \cong [0,1], \quad \sigma(y) = [s]_y \in \mathscr{F}_y \cong [0,1]^2.
\]
Hence,
\[
\textsc{Ps}(\mathrm{Et}(\mathscr{F}))(\{x,y\}) \cong [0,1] \times [0,1],
\]
again with pointwise $\L$ukasiewicz operations.

\medskip
\noindent
Given \( s \in \mathscr{F}(\{x,y\}) \), restriction is defined via:
\[
\begin{aligned}
\textsc{Ps}(\mathrm{Et}(\mathscr{F}))_{\{x,y\} \to \{x\}} &: [0,1] \times [0,1] \to [0,1],\quad (\alpha, \beta) \mapsto \alpha, \\
\textsc{Ps}(\mathrm{Et}(\mathscr{F}))_{\{x,y\} \to \{y\}} &: [0,1] \times [0,1] \to [0,1]^2,\quad \text{identity}.
\end{aligned}
\]

This defines a contravariant functor from \( \mathscr{O}(\mathscr{B}) \) to \( \mathbf{RL} \), the category of residuated lattices, and realizes the presheaf of sections of the \'{e}tal\'{e} space \( \mathrm{Et}(\mathscr{F}) \).

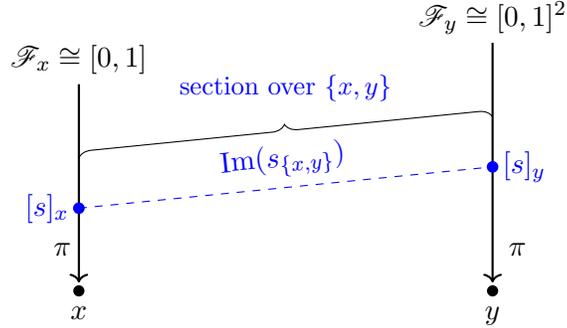
\begin{center}
\begin{tikzpicture}[scale=1.1]
  \node[fill=black, circle, inner sep=1.5pt, label=below:$x$] (x) at (0,0) {};
  \node[fill=black, circle, inner sep=1.5pt, label=below:$y$] (y) at (5,0) {};
  
  \draw[->, thick] (0,2.5) -- (0,0.1);
  \draw[->, thick] (5,3) -- (5,0.1);
  \node at (-0.2,0.5) {$\pi$};
  \node at (5.3,0.5) {$\pi$};

  \node at (0,2.8) {$\mathscr{F}_x \cong [0,1]$};
  \node at (5,3.3) {$\mathscr{F}_y \cong [0,1]^2$};

  \fill[blue] (0,1) circle (2pt) node[left] {$[s]_x$};
  \fill[blue] (5,1.5) circle (2pt) node[right] {$[s]_y$};

  \draw[dashed, blue] (0,1) -- (5,1.5) node[midway, above, sloped] {$\operatorname{Im}(s_{\{x,y\}})$};

  \draw[decorate,decoration={brace,amplitude=5pt},yshift=3pt]
  (0,1.5) -- (5,2) node[midway, above=10pt, blue] {\small section over \(\{x, y\}\)};
\end{tikzpicture}
\end{center}
\end{example}


The following proposition characterizes the elements of the sheafification of a presheaf of residuated lattices.

\begin{proposition}\label{charsectfp}
Let $\mathscr{F}$ be a presheaf of residuated lattices over $\mathscr{B}$, and let $U \in \mathscr{O}(\mathscr{B})$. Then
\[
\mathscr{F}^+(U) = \left\{ \sigma: U \longrightarrow \textsc{Et}(\mathscr{F}) \;\middle|\; 
\begin{array}{l}
\sigma \text{ is a map such that for each } x \in U, \text{ there exists} \\
\text{an open neighborhood } U_x \subseteq U \text{ and } s^x \in \mathscr{F}(U_x) \\
\text{with } \sigma|_{U_x} = s^x_{U_x}
\end{array}
\right\}.
\]
\end{proposition}

\begin{proof}
Let $\Sigma$ denote the set on the right-hand side of the claimed equality. Suppose $\sigma \in \Sigma$. Then for every $x \in U$, there exists an open neighborhood $U_x \subseteq U$ and a section $s^x \in \mathscr{F}(U_x)$ such that $\sigma|_{U_x} = s^x_{U_x}$. Since the family $\{U_x\}_{x \in U}$ is an open cover of $U$, and each restriction $s^x_{U_x}$ is continuous, it follows that $\sigma$ is continuous as well.

Moreover, for each $x \in U$, we have
\[
\pi_{\mathscr{F}}(\sigma(x)) = \pi_{\mathscr{F}}(s^x_{U_x}(x)) = x,
\]
so $\sigma$ is a section of the \'{e}tal\'{e} space $\textsc{Et}(\mathscr{F})$ over $U$.

Conversely, let $\sigma \in \mathscr{F}^+(U)$. For each $x \in U$, the value $\sigma(x)$ belongs to the stalk $\mathscr{F}_x$, so there exists an open neighborhood $W \subseteq U$ of $x$ and a section $s^x \in \mathscr{F}(W)$ such that
\[
\sigma(x) = \mathscr{F}_{W,x}(s^x).
\]
Then $\operatorname{Im}(s^x_W)$ is a neighborhood of $\sigma(x)$ in $\textsc{Et}(\mathscr{F})$, and the restriction $\pi_{\mathscr{F}}|_{\operatorname{Im}(s^x_W)}$ is a homeomorphism. By Lemma~\ref{1shepropo}, there exists an open neighborhood $U_x \subseteq U$ of $x$ such that
\[
\sigma|_{U_x} = \left(\pi_{\mathscr{F}}|_{\operatorname{Im}(s^x_W)}\right)^{-1}\big|_{U_x} = s^x_{U_x}.
\]
Hence, $\sigma \in \Sigma$, completing the proof.
\end{proof}


\section{Relations Between \'Etal\'e Spaces, Presheaves, and Sheaves of Residuated Lattices}\label{sec6}

This section explores the relationships among the categories of étale spaces, presheaves, and sheaves of residuated lattices over a fixed topological space \( \mathscr{B} \). Consider the following diagram that illustrates the core categorical relationships:

\begin{figure}[h!]
\centering
\begin{tikzcd}
\mathbf{RL\text{-}Sheaf}(\mathscr{B}) \arrow[rd, "i" description, hook] &  & \mathbf{RL\text{-}Etale}(\mathscr{B}) \arrow[ll, "\textsc{Ps}" description] \\
 & \mathbf{RL\text{-}PreSheaf}(\mathscr{B}) \arrow[ru, "\textsc{Et}" description] &
\end{tikzcd}
\end{figure}

In the following, we will establish two key results:

\begin{enumerate}
  \item The category \( \mathbf{RL\text{-}Sheaf}(\mathscr{B}) \) is reflective in \( \mathbf{RL\text{-}PreSheaf}(\mathscr{B}) \).
  
  \item The categories \( \mathbf{RL\text{-}Etale}(\mathscr{B}) \) and \( \mathbf{RL\text{-}Sheaf}(\mathscr{B}) \) are equivalent.
\end{enumerate}


Let \( \mathscr{C} \) and \( \mathscr{D} \) be categories. Following \cite[p.~80]{lane1971categories}, an \emph{adjunction} from \( \mathscr{C} \) to \( \mathscr{D} \) is a triple
\[
(\mathscr{L}, \mathscr{R}; \phi): \mathscr{C} \rightharpoonup \mathscr{D},
\]
where \( \mathscr{L}: \mathscr{C} \longrightarrow \mathscr{D} \) and \( \mathscr{R}: \mathscr{D} \longrightarrow \mathscr{C} \) are functors, and \( \phi \) is a natural isomorphism between the bifunctors \( \operatorname{Hom}_{\mathscr{D}}(\mathscr{L}(-), -) \) and \( \operatorname{Hom}_{\mathscr{C}}(-, \mathscr{R}(-)) \). The functor \( \mathscr{L} \) is called a \emph{left adjoint} to \( \mathscr{R} \), and \( \mathscr{R} \) is the \emph{right adjoint} of \( \mathscr{L} \).

\begin{theorem}[\cite{lane1971categories}, \S IV, Theorem 2]\label{equdefadj}
Let \( \mathscr{L}: \mathscr{C} \longrightarrow \mathscr{D} \) and \( \mathscr{R}: \mathscr{D} \longrightarrow \mathscr{C} \) be functors. The following statements are equivalent:
\begin{enumerate}
    \item [$(1)$ \namedlabel{equdefadj1}{$(1)$}] The triple \( (\mathscr{L}, \mathscr{R}; \phi): \mathscr{C} \rightharpoonup \mathscr{D} \) defines an adjunction.
    
    \item [$(2)$ \namedlabel{equdefadj2}{$(2)$}] There exists a natural transformation \( \eta: 1_{\mathscr{C}} \longrightarrow \mathscr{R} \mathscr{L} \) such that for every \( C \in \mathscr{C} \), \( D \in \mathscr{D} \), and morphism \( f: C \to \mathscr{R}(D) \), there exists a unique morphism \( g: \mathscr{L}(C) \to D \) making the following triangle commute:
     \[
      \xymatrix{
\mathscr{L}(C) \ar@{->}[dd]|-{\exists! g} & C \ar@{->}[rr]|-{\eta_{C}} \ar@{->}[rrdd]|-{f} &  & \mathscr{RL}(C) \ar@{->}[dd]|-{\mathscr{R}(g)} \\
 &  &  &  \\
D &  &  & \mathscr{R}(D)
}
\]
\end{enumerate}
\end{theorem}


According to \citet[p.~91]{lane1971categories}, a subcategory \( \mathscr{D} \) of \( \mathscr{C} \) is called \emph{reflective} if the inclusion functor \( i: \mathscr{D} \hookrightarrow \mathscr{C} \) admits a left adjoint \( \mathscr{L}: \mathscr{C} \to \mathscr{D} \). The functor \( \mathscr{L} \) is referred to as a \emph{reflector}, and the adjunction \( (\mathscr{L}, i; \phi): \mathscr{C} \rightharpoonup \mathscr{D} \) is called a \emph{reflection} of \( \mathscr{C} \) into \( \mathscr{D} \).

\begin{theorem}\label{fulsubequdefadj}
Let \( \mathscr{D} \) be a full subcategory of \( \mathscr{C} \). The following are equivalent:
\begin{enumerate}
    \item [$(1)$ \namedlabel{fulsubequdefadj1}{$(1)$}] The subcategory \( \mathscr{D} \) is reflective in \( \mathscr{C} \).
    
    \item [$(2)$ \namedlabel{fulsubequdefadj2}{$(2)$}] For every object \( C \in \mathscr{C} \), there exists an object \( \mathscr{L}(C) \in \mathscr{D} \) and a morphism \( \eta_C: C \to \mathscr{L}(C) \) such that for any morphism \( f: C \to D \) with \( D \in \mathscr{D} \), there exists a unique morphism \( g: \mathscr{L}(C) \to D \) satisfying \( f = g \circ \eta_C \).
\end{enumerate}
\end{theorem}

\begin{proof}
This is an immediate consequence of Theorem~\ref{equdefadj}.
\end{proof}


So far, we have established that \( \mathbf{RL\text{-}Sheaf}(\mathscr{B}) \) forms a full subcategory of \( \mathbf{RL\text{-}PreSheaf}(\mathscr{B}) \) for any given topological space \( \mathscr{B} \). In what follows, we strengthen this result by proving that \( \mathbf{RL\text{-}Sheaf}(\mathscr{B}) \) is in fact a \emph{reflective} subcategory of \( \mathbf{RL\text{-}PreSheaf}(\mathscr{B}) \).


\begin{definition}
Let \( \mathscr{F} \) be a presheaf of residuated lattices over a topological space \( \mathscr{B} \). For each open set \( U \subseteq \mathscr{B} \), define a map
\[
\iota_{\mathscr{F}}(U): \mathscr{F}(U) \longrightarrow \mathscr{F}^{+}(U)
\]
by assigning \( s \mapsto s_U \), where \( s_U: U \to \textsc{Et}(\mathscr{F}) \) is the section defined by \( x \mapsto [s]_x \). 

It can be verified by routine computations that the family \( \iota_{\mathscr{F}} = \{ \iota_{\mathscr{F}}(U) \}_{U \in \mathscr{O}(\mathscr{B})} \) defines a morphism of presheaves of residuated lattices from \( \mathscr{F} \) to its sheafification \( \mathscr{F}^{+} \). This morphism is referred to as the \emph{universal morphism of \( \mathscr{F} \)}.
\end{definition}


\begin{proposition}\label{canomorsheiso}
Let \( \mathscr{F} \) be a sheaf of residuated lattices over a topological space \( \mathscr{B} \). Then the universal morphism \( \iota_{\mathscr{F}}: \mathscr{F} \to \mathscr{F}^+ \) is an isomorphism.
\end{proposition}

\begin{proof}
Let \( U \subseteq \mathscr{B} \) be open. By Proposition~\ref{sectgermequ}, the map \( \iota_{\mathscr{F}}(U): \mathscr{F}(U) \to \mathscr{F}^+(U) \) is injective.

Now let \( \sigma \in \mathscr{F}^+(U) \). By Proposition~\ref{charsectfp}, for every \( x \in U \), there exists an open neighborhood \( U_x \subseteq U \) and a section \( s^x \in \mathscr{F}(U_x) \) such that
\[
\sigma|_{U_x} = (s^x)_{U_x}.
\]
For any \( x, y \in U \), consider the intersection \( U_x \cap U_y \). Let \( z \in U_x \cap U_y \). Then,
\[
\begin{array}{ll}
  r_{U_x \cap U_y, z}(s^x|_{U_x \cap U_y}) & = r_{U_x, z}(s^x) = s^x_{U_x}(z) = \sigma(z) = s^y_{U_y}(z) \\
   & = r_{U_y, z}(s^y) = r_{U_x \cap U_y, z}(s^y|_{U_x \cap U_y}). 
\end{array}
\]

Hence, by Proposition~\ref{sectgermequ}, we conclude that
\[
s^x|_{U_x \cap U_y} = s^y|_{U_x \cap U_y}.
\]
Since the collection \( \{ U_x \}_{x \in U} \) is an open covering of \( U \), and the local sections \( \{ s^x \} \) agree on overlaps, the sheaf condition (condition~\ref{g}) ensures the existence of a unique section \( s \in \mathscr{F}(U) \) such that \( s|_{U_x} = s^x \) for all \( x \in U \). Consequently,
\[
\iota_{\mathscr{F}}(U)(s) = \sigma,
\]
proving that \( \iota_{\mathscr{F}}(U) \) is surjective. Therefore, \( \iota_{\mathscr{F}}(U) \) is a bijection, and \( \iota_{\mathscr{F}} \) is an isomorphism of presheaves.
\end{proof}


\begin{proposition}
Let \( \mathscr{F} \) be a presheaf of residuated lattices over a topological space \( \mathscr{B} \). Then the map
\[
\textsc{Et}(\iota_{\mathscr{F}}): \textsc{Et}(\mathscr{F}) \longrightarrow \textsc{Et}(\mathscr{F}^+)
\]
is an isomorphism of \'{e}tal\'{e} spaces of residuated lattices.
\end{proposition}

\begin{proof}
By Proposition~\ref{etpresheafmor}, the map \( \textsc{Et}(\iota_{\mathscr{F}}) \) is a morphism of \'{e}tal\'{e} spaces of residuated lattices.

To prove injectivity, suppose
\[
\textsc{Et}(\iota_{\mathscr{F}})(\mathscr{F}_{U,x}(s)) = \textsc{Et}(\iota_{\mathscr{F}})(\mathscr{F}_{V,y}(t)).
\]
This means
\[
\mathscr{F}^+_{U,x}(s_U) = \mathscr{F}^+_{V,y}(t_V).
\]
By the definition of germs in the disjoint union, it follows that \( x = y \). Moreover, there exists an open neighborhood \( W \subseteq U \cap V \) of \( x \) such that
\[
\mathscr{F}^+_{UW}(s_U) = \mathscr{F}^+_{VW}(t_V).
\]
Using Proposition~\ref{contavafunctshe}, this implies
\[
s_U(x) = t_V(x),
\]
which means that the original germs are equal:
\[
\mathscr{F}_{U,x}(s) = \mathscr{F}_{V,y}(t).
\]
Hence, \( \textsc{Et}(\iota_{\mathscr{F}}) \) is injective.

To prove surjectivity, let \( \mathscr{F}^+_{U,x}(\sigma) \in \textsc{Et}(\mathscr{F}^+) \). By Proposition~\ref{charsectfp}, there exists an open neighborhood \( U_x \subseteq U \) of \( x \), and a section \( s^x \in \mathscr{F}(U_x) \) such that
\[
\sigma|_{U_x} = (s^x)_{U_x}.
\]
Then,
\[
\textsc{Et}(\iota_{\mathscr{F}})(\mathscr{F}_{U_x,x}(s^x)) = \mathscr{F}^+_{U,x}(\sigma),
\]
showing that \( \textsc{Et}(\iota_{\mathscr{F}}) \) is surjective.

Therefore, \( \textsc{Et}(\iota_{\mathscr{F}}) \) is a bijective morphism of \'{e}tal\'{e} spaces of residuated lattices, and hence an isomorphism.
\end{proof}


\begin{proposition}\label{phipluspro}
Let \( \mathscr{F} \) and \( \mathscr{G} \) be two presheaves of residuated lattices over a topological space \( \mathscr{B} \), and let \( \phi: \mathscr{F} \longrightarrow \mathscr{G} \) be a morphism of presheaves. Then there exists a unique morphism \( \phi^{+}: \mathscr{F}^{+} \longrightarrow \mathscr{G}^{+} \) making the following diagram commute:
\[
\begin{tikzcd}
\mathscr{F} \arrow[rr, "\iota_{\mathscr{F}}" description] \arrow[rdd, "\phi" description] & & \mathscr{F}^{+} \arrow[rdd, "\exists!~\phi^{+}" description,blue] & \\
& & & \\
& \mathscr{G} \arrow[rr, "\iota_{\mathscr{G}}" description] & & \mathscr{G}^{+}
\end{tikzcd}
\]
\end{proposition}

\begin{proof}
Let \( U \subseteq \mathscr{B} \) be open. Define a map
\[
\phi^{+}(U): \mathscr{F}^{+}(U) \longrightarrow \mathscr{G}^{+}(U)
\]
by \( \phi^{+}(U)(\sigma) := \textsc{Et}(\phi) \circ \sigma \), where \( \sigma: U \to \textsc{Et}(\mathscr{F}) \) is a section. Since \( \textsc{Et}(\phi): \textsc{Et}(\mathscr{F}) \to \textsc{Et}(\mathscr{G}) \) is continuous and respects the projection, \( \phi^{+}(U)(\sigma) \) is a section of \( \mathscr{G}^{+} \), and thus \( \phi^{+} \) defines a morphism of presheaves.

Moreover, for each open \( U \subseteq \mathscr{B} \) and \( s \in \mathscr{F}(U) \), we have
\[
\begin{array}{ll}
  \phi^{+}(U)(\iota_{\mathscr{F}}(U)(s)) & = \phi^{+}(U)(s_U) = \textsc{Et}(\phi) \circ s_U \\
   & = (\phi(U)(s))_U = \iota_{\mathscr{G}}(U)(\phi(U)(s)),
\end{array}
\]
so \( \phi^{+} \circ \iota_{\mathscr{F}} = \iota_{\mathscr{G}} \circ \phi \), as desired.

Now suppose \( \psi: \mathscr{F}^{+} \to \mathscr{G}^{+} \) is another morphism of presheaves such that \( \psi \circ \iota_{\mathscr{F}} = \iota_{\mathscr{G}} \circ \phi \). Then applying the \'{e}tal\'{e} functor yields
\[
\textsc{Et}(\psi) \circ \textsc{Et}(\iota_{\mathscr{F}}) = \textsc{Et}(\iota_{\mathscr{G}} \circ \phi) = \textsc{Et}(\phi^{+}) \circ \textsc{Et}(\iota_{\mathscr{F}}).
\]
Since \( \textsc{Et}(\iota_{\mathscr{F}}) \) is an isomorphism by the previous proposition, it follows that \( \textsc{Et}(\psi) = \textsc{Et}(\phi^{+}) \). Because \( \mathscr{G}^{+} \) is a sheaf, Theorem~\ref{preshetalfunc} implies that \( \psi = \phi^{+} \), proving uniqueness.
\end{proof}


The following theorem establishes that the category of sheaves of residuated lattices over a topological space is reflective in the corresponding category of presheaves. Specifically, it shows that the composite functor \( \textsc{Ps} \circ \textsc{Et} \) serves as a left adjoint to the inclusion functor.

\begin{theorem}\label{reflshbpshb}
For a topological space $\mathscr{B}$, \( \mathbf{RL\text{-}Sheaf}(\mathscr{B}) \) is reflective in \( \mathbf{RL\text{-}PreSheaf}(\mathscr{B}) \).
\end{theorem}

\begin{proof}
Let \( \mathscr{F} \) be a presheaf and \( \mathscr{G} \) a sheaf of residuated lattices over \( \mathscr{B} \), and suppose \( \phi: \mathscr{F} \to \mathscr{G} \) is a morphism of presheaves.

By Propositions~\ref{canomorsheiso} and~\ref{phipluspro}, the universal morphism \( \iota_{\mathscr{G}}: \mathscr{G} \to \mathscr{G}^+ \) is an isomorphism, and there exists a unique morphism \( \phi^+: \mathscr{F}^+ \to \mathscr{G}^+ \) such that
\[
\phi^+ \circ \iota_{\mathscr{F}} = \iota_{\mathscr{G}} \circ \phi.
\]
Since \( \iota_{\mathscr{G}} \) is an isomorphism, we can define
\[
\hat{\phi} := \iota_{\mathscr{G}}^{-1} \circ \phi^+ : \mathscr{F}^+ \longrightarrow \mathscr{G},
\]
and this morphism satisfies \( \hat{\phi} \circ \iota_{\mathscr{F}} = \phi \), as shown in the following diagram:

\[
\begin{tikzcd}
\mathscr{F} \arrow[rr, "\iota_{\mathscr{F}}" description] \arrow[rdd, "\phi" description] & & \mathscr{F}^{+} \arrow[rdd, "\phi^{+}" description] \arrow[ldd, "\hat{\phi} = \iota_{\mathscr{G}}^{-1} \circ \phi^{+}" description, dashed,blue] & \\
& & & \\
& \mathscr{G} \arrow[rr, "\iota_{\mathscr{G}}" description] & & \mathscr{G}^{+}
\end{tikzcd}
\]

Finally, the uniqueness of \( \hat{\phi} \) follows from the uniqueness of \( \phi^+ \), as established in Proposition~\ref{phipluspro}. Hence, the functor \( \mathscr{F} \mapsto \mathscr{F}^+ \) is left adjoint to the inclusion \( i \), proving the reflectivity.
\end{proof}

\begin{remark}
Theorem~\ref{reflshbpshb} highlights a fundamental principle: every presheaf of residuated lattices can be canonically completed into a sheaf using a geometrically and categorically motivated process. Even if a presheaf contains inconsistencies or data that fails to glue, the sheafification functor
\[
\mathscr{F} \mapsto \mathscr{F}^+ := \textsc{Ps} \circ \textsc{Et}(\mathscr{F})
\]
acts as a \emph{repair mechanism}, transforming this data into a globally coherent, semantically faithful sheaf.

The construction first passes through the \'{e}tal\'{e} space \( \mathrm{Et}(\mathscr{F}) \), offering a geometric view of how local behaviors cluster around points. The functor \( \textsc{Ps} \) then extracts continuous sections from this space, enforcing the gluing conditions that define a sheaf. The result is a sheaf \( \mathscr{F}^+ \) that best approximates \( \mathscr{F} \) while preserving its algebraic structure.

\begin{center}
\begin{tikzpicture}[node distance=3.2cm, auto]
  \node (F) {\( \mathscr{F} \)};
  \node (EtF) [right of=F] {\( \mathrm{Et}(\mathscr{F}) \)};
  \node (Fplus) [right of=EtF] {\( \mathscr{F}^+ \)};
  \node (G) [below of=Fplus, node distance=2.3cm] {\( \mathscr{G} \)};
  
  \draw[->, thick] (F) to node[above] {\scriptsize \( \textsc{Et} \)} (EtF);
  \draw[->, thick] (EtF) to node[above] {\scriptsize \( \textsc{Ps} \)} (Fplus);
  \draw[->, thick, dashed,blue] (Fplus) to node[right] {\scriptsize \( \hat{\phi} \)} (G);
  \draw[->, thick] (F) to[out=-30, in=180] node[below] {\scriptsize \( \phi \)} (G);
\end{tikzpicture}
\end{center}

This diagram illustrates the universal factorization property: any morphism \( \phi: \mathscr{F} \to \mathscr{G} \) into a sheaf uniquely factors through \( \mathscr{F}^+ \), making the category \( \mathbf{RL\text{-}Sheaf}(\mathscr{B}) \) reflective in \( \mathbf{RL\text{-}PreSheaf}(\mathscr{B}) \).

From a logical and semantic perspective, this tells us that sheaves are not restrictive—they are the \emph{best possible upgrade} of presheaves. The functor \( \mathscr{F} \mapsto \mathscr{F}^+ \) preserves all meaningful local information while discarding only what cannot be made coherent globally. In substructural and fuzzy logics, this guarantees that locally defined truth values or algebraic structures can always be extended to global ones in a principled, universal manner.

In summary, the theorem assures us that the transition from local to global is not only possible but optimally structured. Sheafification offers a canonical way to pass from potentially inconsistent local data to a globally consistent and semantically meaningful whole—essential for applications in topology, algebra, and logic.
\end{remark}


Let $\mathscr{C}$ and $\mathscr{D}$ be categories. Recall \cite[p.93]{lane1971categories} that a functor $S:\mathscr{C} \longrightarrow \mathscr{D}$ is an \textit{equivalence of categories} (and that $\mathscr{C}$ and $\mathscr{D}$ are \textit{equivalent}) if there exists a functor $T:\mathscr{D} \longrightarrow \mathscr{C}$ and natural isomorphisms $TS \cong 1_{\mathscr{C}} : \mathscr{C} \longrightarrow \mathscr{C}$ and $ST \cong 1_{\mathscr{D}} : \mathscr{D} \longrightarrow \mathscr{D}$. In this case, the functor $T$ is also an equivalence of categories.

We now show that the category of \'{e}tal\'{e} spaces of residuated lattices over a given topological space is equivalent to the category of sheaves of residuated lattices over the same space. For a fixed space $\mathscr{B}$, we prove that the functor $\textsc{Ps}$ is an equivalence with quasi-inverse given by $\textsc{Et} \circ i$.

\begin{theorem}\label{eqetshtheo}
Let $\mathscr{B}$ be a topological space. Then the categories $\mathbf{RL\text{-}Etale}(\mathscr{B})$ and $\mathbf{RL\text{-}Sheaf}(\mathscr{B})$ are equivalent.
\end{theorem}

\begin{proof}
Let $\mathscr{T} \overset{\pi}{\downarrow} \mathscr{B}$ be an \'{e}tal\'{e} space of residuated lattices. Define the map 
\[
\varepsilon_{\mathscr{T}} : \textsc{Et}(\textsc{Ps}(\mathscr{T})) \longrightarrow \mathscr{T}
\]
by 
\[
\textsc{Ps}(\mathscr{T})_{U,x}(\sigma) \mapsto \sigma(x).
\]
Suppose $\textsc{Ps}(\mathscr{T})_{U,x}(\sigma) = \textsc{Ps}(\mathscr{T})_{V,x}(\rho)$. Then there exists a neighborhood $W \subseteq U \cap V$ of $x$ such that $\textsc{Ps}(\mathscr{T})_{UW}(\sigma) = \textsc{Ps}(\mathscr{T})_{VW}(\rho)$, which implies $\sigma(x) = \rho(x)$.

Let $t \in \mathscr{T}$. By Lemma \ref{shepropo}, there exists an open neighborhood $U$ of $\pi(t)$ and a section $\sigma$ over $U$ such that $\sigma(\pi(t)) = t$. Then
\[
\varepsilon_{\mathscr{T}}(\textsc{Ps}(\mathscr{T})_{U,\pi(t)}(\sigma)) = t,
\]
showing that $\varepsilon_{\mathscr{T}}$ is surjective.

Next, assume 
\[
\varepsilon_{\mathscr{T}}(\textsc{Ps}(\mathscr{T})_{U,x}(\sigma)) = \varepsilon_{\mathscr{T}}(\textsc{Ps}(\mathscr{T})_{V,y}(\rho)).
\]
Then $\sigma(x) = \rho(y)$, and consequently $x = \pi(\sigma(x)) = \pi(\rho(y)) = y$. By Proposition \ref{sheprop}\ref{sheprop1}, the equalizer $W = \textsc{Eq}(\sigma, \rho)$ is an open set contained in both $U$ and $V$, and 
\[
\textsc{Ps}(\mathscr{T})_{UW}(\sigma) = \textsc{Ps}(\mathscr{T})_{VW}(\rho),
\]
so 
\[
\textsc{Ps}(\mathscr{T})_{U,x}(\sigma) = \textsc{Ps}(\mathscr{T})_{V,x}(\rho).
\]
Thus, $\varepsilon_{\mathscr{T}}$ is injective, and hence a bijection.

Moreover, for each $b \in \mathscr{B}$, we have
\[
\varepsilon_{\mathscr{T}}(\textsc{Et}(\textsc{Ps}(\mathscr{T}))_{b}) \subseteq \mathscr{T}_{b}.
\]
Let $U$ be an open set in $\mathscr{B}$ and let $\sigma \in \textsc{Ps}(\mathscr{T})(U) = \Gamma(U, \mathscr{T})$. Then
\[
\varepsilon_{\mathscr{T}}(\mathrm{Im}(\sigma_U)) = \sigma(U).
\]
By Propositions \ref{conopnloc} and \ref{sheprop}\ref{sheprop3}, it follows that $\varepsilon_{\mathscr{T}}$ is an isomorphism of \'{e}tal\'{e} spaces.

Now fix $x \in \mathscr{B}$. The following sequence shows that $\varepsilon_{\mathscr{T},x}$ preserves the join operation $\vee$:
\[
\begin{array}{ll}
\varepsilon_{\mathscr{T},x}(\textsc{Ps}(\mathscr{T})_{U,x}(\sigma) \vee \textsc{Ps}(\mathscr{T})_{V,x}(\tau)) 
& = \varepsilon_{\mathscr{T},x}(\textsc{Ps}(\mathscr{T})_{U \cap V, x}(\sigma|_{U \cap V} \vee \tau|_{U \cap V})) \\
& = (\sigma|_{U \cap V} \vee \tau|_{U \cap V})(x) \\
& = \sigma(x) \vee \tau(x) \\
& = \varepsilon_{\mathscr{T},x}(\textsc{Ps}(\mathscr{T})_{U,x}(\sigma)) \vee \varepsilon_{\mathscr{T},x}(\textsc{Ps}(\mathscr{T})_{V,x}(\tau)).
\end{array}
\]
A similar argument shows that $\varepsilon_{\mathscr{T},x}$ preserves $\wedge$, $\odot$, and $\rightarrow$. Hence, $\varepsilon_{\mathscr{T}}$ is an isomorphism of \'{e}tal\'{e} spaces of residuated lattices.

Now let $\mathscr{T} \overset{\pi}{\downarrow} \mathscr{B}$ and $\mathscr{S} \overset{\phi}{\downarrow} \mathscr{B}$ be two \'{e}tal\'{e} spaces of residuated lattices, and let $h : \mathscr{T} \rightarrow \mathscr{S}$ be a morphism of such spaces. Then we have:
\[
\begin{array}{ll}
h \circ \varepsilon_{\mathscr{T}}(\textsc{Ps}(\mathscr{T})_{U,x}(\sigma)) 
& = h(\sigma(x)) \\
& = \varepsilon_{\mathscr{S}}(\textsc{Ps}(\mathscr{S})_{U,x}(h \circ \sigma)) \\
& = \varepsilon_{\mathscr{S}}(\textsc{Et}(\textsc{Ps}(h))(\textsc{Ps}(\mathscr{T})_{U,x}(\sigma))).
\end{array}
\]
Therefore, the following diagram commutes:
\[
\begin{tikzcd}
\textsc{Et}(\textsc{Ps}(\mathscr{T})) \arrow[dd, "\varepsilon_{\mathscr{T}}" description] \arrow[rrrr, "\textsc{Et}(\textsc{Ps}(h))" description,blue] & & & & \textsc{Et}(\textsc{Ps}(\mathscr{S})) \arrow[dd, "\varepsilon_{\mathscr{S}}" description] \\
& & & & \\
\mathscr{T} \arrow[rrrr, "h" description] & & & & \mathscr{S}
\end{tikzcd}
\]
Thus, $(\textsc{Et} \circ i) \circ \textsc{Ps} \cong 1_{\textsc{Et}(\mathscr{B})}$.

Now let $(\mathscr{F}, \mathscr{B})$ be a sheaf of residuated lattices. By Proposition \ref{canomorsheiso}, the canonical morphism 
\[
l_{\mathscr{F}} : \mathscr{F} \longrightarrow \textsc{Ps}(\textsc{Et}(\mathscr{F}))
\]
is an isomorphism. Moreover, Proposition \ref{phipluspro} ensures the coherence of morphisms under composition. Hence,
\[
\textsc{Ps} \circ (\textsc{Et} \circ i) \cong 1_{\mathbf{RL\text{-}Sheaf}(\mathscr{B})}.
\]
\end{proof}


Theorem~\ref{eqetshtheo} shows that the category of \'{e}tal\'{e} spaces of residuated lattices over a topological space \( \mathscr{B} \) is equivalent to the category of sheaves of residuated lattices over the same base. This result is conceptually important because it gives us two equivalent but distinct ways to think about the same mathematical objects—one geometric, the other algebraic.

From the \'{e}tal\'{e} space perspective, a sheaf is visualized as a bundle of stalks (local structures) topologically varying over the base space. From the sheaf-theoretic perspective, we think in terms of assigning consistent local data (sections) over open sets, subject to gluing conditions. The equivalence says that these two pictures are not just similar—they are categorically the same.

\begin{center}
\begin{tikzcd}[column sep=large, row sep=large]
\mathbf{RL\text{-}Etale}(\mathscr{B}) \arrow[r, shift left=1.2ex, "\textsc{Ps}"] 
& \mathbf{RL\text{-}Sheaf}(\mathscr{B}) \arrow[l, shift left=1.2ex, "\textsc{Et} \circ i"]
\end{tikzcd}
\end{center}

This adjoint equivalence means that:
- Every \'{e}tal\'{e} space can be converted into a sheaf of sections;
- Every sheaf arises (up to isomorphism) from such a space of stalks;
- And natural isomorphisms \( \textsc{Ps} \circ (\textsc{Et} \circ i) \cong 1 \) and \( (\textsc{Et} \circ i) \circ \textsc{Ps} \cong 1 \) hold.

Practically, this means we can freely switch between a local, spatial model of a sheaf (as an \'{e}tal\'{e} space) and an algebraic model (as a functor satisfying sheaf conditions). The functor \( \textsc{Ps} \) extracts the sheaf of continuous sections from a geometric \'{e}tal\'{e} space, while \( \textsc{Et} \) reconstructs the space of germs from a sheaf. The natural isomorphisms in both directions guarantee that no information is lost in translation.

This equivalence is especially useful in logical and semantic applications, where one often builds models algebraically (as presheaves or sheaves) but wants to interpret them geometrically (as variable structures over a space). Étale spaces offer an intuitive framework for understanding how local truth-values or logical formulas vary across space or context, while sheaves ensure that those variations are coherent. In this sense, the theorem provides a formal bridge between algebraic reasoning and spatial semantics.


\section{Conclusion and Future Work}\label{sec7}

This section presents the conclusion of the article, offers a logical interpretation of the main definitions and theorems, and outlines potential directions for future research inspired by the results of this work.

\subsection*{Conclusion}

In this work, we introduced and developed the theory of \emph{\'{e}tal\'{e} spaces of residuated lattices}, expanding the interaction between algebraic logic and topology. Building upon classical concepts from sheaf theory, we defined \'{e}tal\'{e} spaces where each stalk is a residuated lattice and where local sections inherit algebraic structure via pointwise operations. We demonstrated that such spaces form a well-defined subcategory of the category of \'{e}tal\'{e} spaces over a fixed topological base space $B$, denoted as $\mathrm{RL\text{-}Etale}(B)$.

Moreover, we explored the categorical underpinnings of this construction. Morphisms between \'{e}tal\'{e} spaces of residuated lattices were carefully defined to preserve the algebraic structure at each stalk, ensuring the compositional and identity properties necessary for a subcategory. We also provided concrete examples to illustrate how classical residuated lattices (such as $A_4$, $A_6$, and $A_8$) yield interesting topological and categorical structures.

In the second part of the paper, we introduced the notion of \emph{presheaves and sheaves of residuated lattices}, framed as contravariant functors from the lattice of open sets $\mathcal{O}(B)$ to the category $\mathrm{RL}$ of residuated lattices. We established the connection between \'{e}tal\'{e} spaces and sheaves through a detailed analysis of the gluing and separation properties. These results provide a logical and categorical foundation for modeling local-to-global reasoning in many-valued and substructural logics.

\subsection*{Logical Interpretation}

From a logical and semantic viewpoint, the constructions in this paper reflect a topological semantics for non-classical logics, particularly those governed by residuated lattices, such as substructural logics and fuzzy logics. In this framework, each point in the base topological space corresponds to a local logical environment, and the stalk over that point provides the algebraic structure of truth values valid in that context.

The sections of an \'{e}tal\'{e} space represent contextually valid assignments or formulas, and the continuity conditions of the \'{e}tal\'{e} projection ensure that these assignments vary coherently over open neighborhoods. The logical connectives (such as conjunction $\wedge$, implication $\to$, and fusion $\odot$) are interpreted as continuous operations defined fiberwise, making the entire structure a semantic domain for interpreting context-sensitive logical formulas.

Presheaves and sheaves refine this perspective further: presheaves allow for modeling partially defined logical information across open sets, while sheaves enforce a consistency condition, ensuring that locally compatible data can be uniquely amalgamated. The equivalence between certain categories of sheaves and \'{e}tal\'{e} spaces thus provides a powerful geometric realization of logical systems, aligning with categorical logic and topos theory.

\subsection*{Future Work}

Several directions remain open for future exploration:

\begin{itemize}
    \item \textbf{Non-commutative generalizations:} A natural next step is to generalize the current framework to \emph{non-commutative residuated lattices}. This would enable the modeling of richer logical systems, such as relevance logics and non-commutative linear logics, which require relaxing the commutativity of fusion $\odot$. Extending \'{e}tal\'{e} space constructions and sheaf semantics to this setting will likely require new categorical tools, such as enriched categories or bicategories.

    \item \textbf{Internal sheaves and logical toposes:} Another line of inquiry involves investigating internal sheaves of residuated lattices within suitable logical or algebraic toposes. This could allow for internalizing our constructions and studying their behavior under logical functors, with potential applications in categorical logic and model theory.

    \item \textbf{Dynamical and epistemic systems:} Sheaf-theoretic methods are increasingly used in the modeling of distributed and dynamic systems, particularly in multi-agent and temporal logics. Our framework can be extended to incorporate temporal evolution or agent-indexed residuated lattices, potentially giving rise to sheaves on fibered spaces or indexed categories.

    \item \textbf{Logical dualities and completeness:} Finally, it is worth investigating whether the constructions developed here admit duality theorems (in the spirit of Stone or Esakia duality), and whether the sheaf representations introduced can yield completeness results for certain fragments of substructural logic. This may involve constructing canonical sheaves over suitable spectral spaces associated with the logical syntax.

\end{itemize}

Overall, this work lays the foundation for a broader research program at the intersection of logic, topology, and algebra, aimed at deepening our understanding of context-dependent truth and local reasoning in logical systems.


\end{document}